\documentclass{article}
\usepackage{caption}
\usepackage{geometry}
\usepackage{graphicx}

\usepackage[T1]{fontenc}
\usepackage{amssymb}
\usepackage{amsmath,cases}
\usepackage{amsthm}
\usepackage{floatrow}
\usepackage{mathtools}
\usepackage{enumitem}
\usepackage{xcolor}
\usepackage{textcomp}
\usepackage{mathrsfs}
\usepackage{bbm}
\usepackage{bm}
\usepackage{chemarrow}
\usepackage{stmaryrd}
\usepackage{url}
\usepackage{calc}
\usepackage{ifthen}
\usepackage{psfrag}
\usepackage[active]{srcltx}   
\usepackage[all,cmtip]{xy}
\usepackage{pdfsync}
\usepackage[refpage]{nomencl}  
\usepackage{yhmath}
\usepackage{framed}
\usepackage{todonotes}

\newcommand{\eps}{\varepsilon}

\newcommand{\de}{\mathmbox{\mathrel{\mathop:}\hspace*{-.6pt}=}}
\newcommand{\ooo}{\bullet}
\newcommand{\oov}{{\mathrlap{\textcolor{vert}{\ooo}}{\circ}}}
\newcommand{\sqr}{{\mathrlap{\scalebox{.37}{\rotatebox{45}{\textcolor{rouge}{$\blacksquare$}}}}{\scalebox{.37}{\rotatebox{45}{$\square$}}}}}
\newcommand{\si}{\sigma}
\newcommand{\un}[1]{\mathbbm{1}_{#1}}

\newcommand{\va}{\vec{a}}
\newcommand{\vb}{\vec{\smash{b}}}
\newcommand{\vc}{\vec{c}}
\newcommand{\vrho}{\vec{\rho}}
\newcommand{\vtau}{\vec{\tau}}

\newcommand{\N}{\mathbb{N}}
\newcommand{\Z}{\mathbb{Z}}

\newcommand{\cB}{\mathcal{B}}
\newcommand{\cC}{\mathcal{C}}
\newcommand{\EE}{\mathcal{E}}
\newcommand{\calG}{\mathcal{G}}
\newcommand{\M}{\mathcal{M}}
\newcommand{\cO}{\mathcal{O}}
\newcommand{\OO}{\mathcal{O}}
\newcommand{\cS}{\mathcal{S}}
\newcommand{\T}{\mathcal{T}}
\newcommand{\cU}{\mathcal{U}}

\newcommand{\lab}{\mathfrak{l}}
\newcommand{\m}{\mathfrak{m}}
\newcommand{\fM}{\mathfrak{M}}
\newcommand{\q}{\mathfrak{q}}
\newcommand{\s}{\mathfrak{s}}
\newcommand{\Sg}{\mathfrak{S}}
\newcommand{\fu}{\mathfrak{u}}

\DeclareMathOperator{\suc}{succ}
\DeclareMathOperator{\pr}{par}
\DeclareMathOperator{\opp}{opp}
\DeclareMathOperator{\upd}{upd}
\DeclareMathOperator{\update}{update}
\newcommand{\im}{\operatorname{ind}^{\text-}}
\newcommand{\ip}{\operatorname{ind}^{\text+}}

\definecolor{vert}{RGB}{85,170,0}
\definecolor{rouge}{RGB}{255,170,170}

\theoremstyle{plain}
\newtheorem{thm}{Theorem}
\newtheorem{lem}[thm]{Lemma}
\newtheorem{prop}[thm]{Proposition}

\newtheorem{corol}[thm]{Corollary}
\newtheorem{defi}{Definition}

\theoremstyle{definition}
\newtheorem*{rem}{Remark}

\newtheorem*{note}{Note}
\newtheorem*{ack}{Acknowledgment}

\newcommand{\edg}[2]{
\begin{tikzpicture}[baseline=(X.base), inner sep=0, outer sep=0]
	\node[draw,circle,minimum height=4.5mm] (X) at (0,0) {\footnotesize #1};
	\node[draw,circle,minimum height=4.5mm] (Y) at (.8,0) {\footnotesize #2};
	\draw (X)--(Y);
  \end{tikzpicture}
}

\newcommand{\corn}[3]{
\begin{tikzpicture}[baseline=(X.base), inner sep=0, outer sep=0]
	\node[draw,circle,minimum height=4.5mm] (X) at (0,0) {\footnotesize #1};
	\node[draw,circle,minimum height=4.5mm] (Y) at (.8,0) {\footnotesize #2};
	\node[draw,circle,minimum height=4.5mm] (Z) at (1.6,0) {\footnotesize #3};
	\draw (X)--(Y)--(Z);
  \end{tikzpicture}
}

\newcounter{step}
\makenomenclature

\makeatletter
\newcommand \Dotfill {\leavevmode \leaders \hb@xt@ .72em{\hss .\hss } \hfill \kern \z@}
\makeatother
\renewcommand*{\pagedeclaration}[1]{\Dotfill\ifthenelse{#1<10}{\hspace{.72em}}{\hspace{.2em}}\hyperpage{#1}}
\usepackage[colorlinks=true,bookmarksnumbered=true, pdfauthor={J\'er\'emie Bettinelli}]{hyperref}
\newenvironment{pre}[1][\proofname]{%
  \proof[#1]%
}{\endproof}

\makeatletter
\renewcommand*{\@fnsymbol}[1]{\ensuremath{\ifcase#1\or \dagger\or \ddagger\or
   \mathsection\or \mathparagraph\or \|\or **\or \dagger\dagger
   \or \ddagger\ddagger \else\@ctrerr\fi}}
\makeatother
\title{A bijection for nonorientable general maps}
\author{J\'er\'emie \textsc{Bettinelli}\thanks{LIX, cnrs, \'Ecole polytechnique, Institut Polytechnique de Paris.}}

\begin{document}
\maketitle

\begin{abstract}
We give a different presentation of a recent bijection due to Chapuy and Do{\l}{\k{e}}ga for nonorientable bipartite quadrangulations and we extend it to the case of nonorientable general maps. This can be seen as a Bouttier--Di Francesco--Guitter-like generalization of the Cori--Vauquelin--Schaeffer bijection in the context of general nonorientable surfaces. In the particular case of triangulations, the encoding objects take a particularly simple form and this allows us to recover a famous asymptotic enumeration formula found by Gao.
 
\end{abstract}


\section{Introduction}

\subsection{Motivation}

This paper is an improved extended version of~\cite{bettinelli15fpsac}. The study of maps has seen tremendous developments in the past few decades. One of the reasons is that they provide natural discrete versions of a given surface. In particular, when taken according to a well-chosen natural probability distribution, it has been shown for several models that a random map converges (after scaling, in a certain sense) toward a limiting object. This limiting object is a random metric space and has (almost surely) the same topology as the surface on which the considered maps are drawn. It is called the \emph{Brownian sphere} or \emph{Brownian map} when the surface is the sphere, and the \emph{Brownian $\cS$} for a general orientable surface~$\cS$.

In the most-studied case of the sphere, it has been shown \cite{legall11ubm,miermont11bms} for example that a uniform quadrangulation (map with all faces of degree~$4$) with~$n$ faces converges to the Brownian sphere as $n\to\infty$. See also \cite{beltranlg,addarioalbenque2013simple,BeJaMi14,abr13,Mar16,addarioalbenque2020} for other natural models of maps drawn on the sphere that exhibit the same behavior. An important aspect of these results is that the limiting object is universal, in the sense that it is always the Brownian sphere (up to a scale constant), independently of the model one considers.

The case of a more general compact orientable surface (with a boundary allowed) has been studied, mostly in the context of uniform quadrangulations: partial convergence has been established in a series of papers ending with~\cite{bettinelli14gbs}, and full convergence is under investigation~\cite{BeMi15II}. The full convergence in the particular case of the disk has recently been shown in~\cite{BeMi15I}, where many more models are also considered.

All the previously mentioned results strongly rely on powerful bijective encodings of the considered maps. It turns out that quadrangulations are particularly well behaved with respect to these bijective encodings and this is the main reason why they are usually the first to be studied. However, if one wishes to study other models and, in particular, surfaces with a boundary, one needs more general bijective encodings. In the case of compact orientable surfaces, the so-called \emph{Schaeffer-like bijections} \cite{cori81planar,schaeffer98cac,bouttier04pml,chapuy07brm,miermont09trm,ambjornbudd} allow one to conduct most studies. Note that, in certain cases, bijections of a different kind (\emph{blossoming bijections})~\cite{schaeffer97bij,poulscha06,FuPoSc09,AlPo15,Le19,DoLe20} have also been used~\cite{addarioalbenque2013simple,addarioalbenque2020}.

\begin{figure}[ht!]
		\psfrag{v}[][]{$v^\ooo$}
		\psfrag{0}[B][B]{}
		\psfrag{1}[B][B]{}
		\psfrag{2}[B][B]{}
		\psfrag{3}[B][B]{}
		\psfrag{4}[B][B]{}
	\centering\includegraphics[width=8cm]{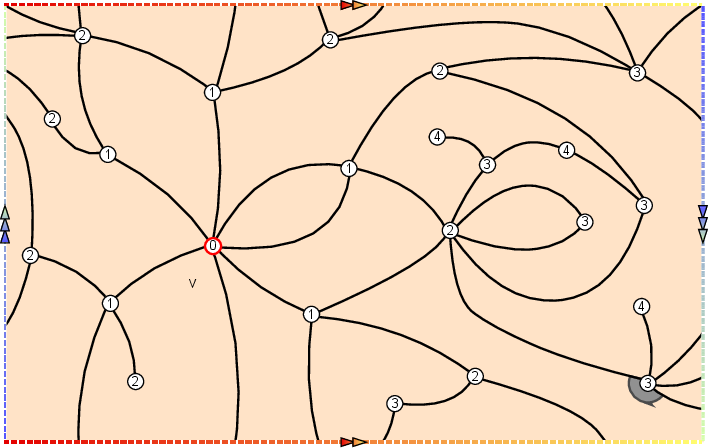}
	\caption{A pointed bipartite map of the Klein bottle with $16$ faces. The root (see below) is the oriented corner symbolized by the gray half arrowhead.}
	\label{exmap}
\end{figure}

Until very recently, no such bijections were known in the case of a nonorientable surface. In~\cite{ChDo15bij}, Chapuy and Do{\l}{\k{e}}ga took the first step by exhibiting a bijection allowing to encode nonorientable bipartite quadrangulations. In this work, we give a bijection inspired by their construction, which provides an explicit construction for pointed nonorientable general maps\footnote{A \emph{pointed map} is a map given with a distinguished vertex. From a combinatorial point of view, it might not seem to make much difference, but the bijections for pointed maps turn out to be better behaved for probabilistic applications.}. We use a rather different presentation than in~\cite{ChDo15bij}, which will be more suited for our generalization: in particular, we introduce the notion of \emph{level loop}, which replaces the notion of dual exploration graph used in~\cite{ChDo15bij}. These works lay the bases for the future study of nonorientable Brownian surfaces~\cite{BeChDo15}. Figure~\ref{exmap} shows an example of a pointed bipartite map of the Klein bottle.

\bigskip
Another cause of interest for maps is their remarkable enumerative properties. In fact, although maps are intricate objects by nature, many classes of them possess a quite simple enumerative structure. Thanks to different involved enumeration techniques (generating functions, matrix integrals, algebraic combinatorics), many classes of maps have been enumerated and map enumeration has become over the years a full-fledged research domain.

In the case of the sphere, Tutte~\cite{tutte63cpm} gave a very simple closed formula for the number of rooted maps with a given number of edges. A bijective proof of this formula was given by Cori and Vauquelin~\cite{cori81planar} and later popularized by Schaeffer~\cite{schaeffer98cac}. It relies on their so-called \emph{Cori--Vauquelin--Schaeffer bijection} encoding quadrangulations of the sphere with trees whose vertices carry integer labels satisfying local constraints. For more general surfaces, Bender and Canfield~\cite{bender86anr} showed that the number of rooted maps with~$n$ edges on a given surface (orientable or not) is asymptotically equal to a constant times $n^{5(h-1)/2} 12^n$, where $h$ is the type of the considered surface and the constant depends on the surface. Extending the Cori--Vauquelin--Schaeffer bijection, a combinatorial interpretation of this fact in the orientable case was given by Chapuy, Marcus and Schaeffer~\cite{chapuy07brm}. Their approach relies on a bijection between bipartite quadrangulations (it is a classical simple fact that bipartite quadrangulations are in bijection with general maps) and one-face maps of the same surface, whose vertices carry integer labels satisfying some local constraints.

In parallel, Bouttier, Di Francesco and Guitter~\cite{bouttier04pml} extended the original Cori--Vauquelin--Schaeffer bijection to encode maps of the sphere with an arbitrary face degree distribution. Unifying both aforementioned extensions, Chapuy~\cite{chapuy09ori} proved similar asymptotic enumeration results for more families of maps of an orientable surface.

Nonorientability does not cause too much difficulties for generating function approaches. In addition to Bender and Canfield's results, we may cite the work of Gao, who showed~\cite{gao91} that the number of rooted triangulations with~$n$ edges on a given surface (orientable or not) is asymptotically equal to a constant (depending on the surface) times $n^{5(h-1)/2} (12\sqrt{3})^n$. He also studied the algebraicity of the generating function 
of rooted maps of a given surface with face degree constraints~\cite{gao93number}.

In their recent work~\cite{ChDo15bij}, Chapuy and Do{\l}{\k{e}}ga extended the construction of~\cite{chapuy07brm} to bipartite nonorientable quadrangulations. In this paper, we give a different construction of their bijection and extend it by an approach reminiscent of~\cite{bouttier04pml}. In order to achieve this goal, we somehow fix a local orientation of the surface via a global process, in the sense that the process uses the information of the whole map. As a result, the constraints satisfied by the labels of the encoding objects are also \emph{global} and this prevents us from giving a simple characterization of these objects. In the very particular case of bipartite quadrangulations, the global constraints can be expressed as local constraints and the encoding objects (so-called \emph{well-labeled unicellular maps}) take a simple form. This allowed Chapuy and Do{\l}{\k{e}}ga~\cite{ChDo15bij} to give a combinatorial interpretation of Bender and Canfield's asymptotic formula. In the particular case of triangulations, the same miracle occurs and we are able to give a combinatorial interpretation of the results of~\cite{gao91}.

\subsection{First definitions}

In the present paper, we work on a compact surface without boundary. Recall that, by the classification theorem, it is either orientable and homeomorphic to the surface obtained by adding~$h$ handles to the sphere for some $h\in\{0,1,\ldots\}$ (sphere, torus, double torus, etc.), or nonorientable and homeomorphic to the surface obtained by adding~$2h$ cross-caps to the sphere for some $h\in\{1/2,1,3/2,\ldots\}$ (projective plane, Klein bottle, etc.). The number~$h$ is called the \emph{type} of the surface.

\begin{framed}
From now on and until the end of the paper, we fix such a surface~$\cS$, orientable or not.
\end{framed}
\nomenclature[1a]{$\cS$}{compact surface without boundary on which we work, fixed once and for all}%
\nomenclature[1b]{$h$}{type of $\cS$}%

\begin{figure}[ht!]
		\psfrag{r}[][][.8]{$\vrho$}
		\psfrag{s}[][][.8]{$\vec\gamma$}
		\psfrag{v}[][][.8]{$v$}
		\psfrag{e}[][][.8]{$e$}
		\psfrag{c}[][][.8]{$r$}
		\psfrag{a}[][][.7]{$\vc$}
		\psfrag{p}[][][.7]{$\varphi(\vc)$}
		\psfrag{d}[][][.7]{$\sigma(\vc)$}
	\centering\includegraphics[width=12cm]{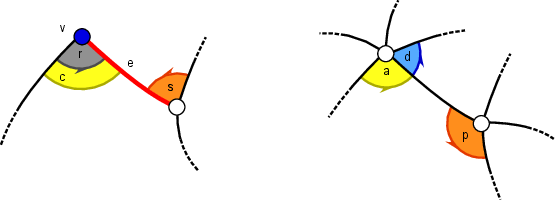}
	\caption{\textbf{Left.} The (gray) root~$\vrho$, the (yellow) root corner~$r$, the (blue) root vertex~$v$ and the (red) root edge~$e$. The root of the root flipped map is~$\vec\gamma$ (in orange). \textbf{Right.} An oriented corner and its images by~$\sigma$ and~$\varphi$.}
	\label{rootthgs}
\end{figure}

A \emph{map} is a cellular embedding of a finite graph (possibly with multiple edges and loops) into~$\cS$, considered up to homeomorphisms. Cellular means that the connected components of the complement of edges, called \emph{faces}, are homeomorphic to 2-dimensional open disks. Note that, although~$\cS$ might not be orientable, each face is orientable, so it will make sense to follow the border of a face, provided a starting orientation is given. A \emph{corner} is an angular sector determined by two consecutive half-edges incident to the same vertex and to the same face. The \emph{degree} of a face is its number of corners. By convention, all the maps we consider are \emph{rooted}, that is, given with a distinguished \textbf{oriented} corner called the \emph{root}, usually denoted by~$\vrho$ (for a general map) or~$\vtau$ (for a unicellular map). It will be represented by a small angular sector with an arrow on the pictures. The (nonoriented) corner corresponding to the root will be called the \emph{root corner}, the vertex incident to the root will be called the \emph{root vertex} and the edge incident to the root will be called the \emph{root edge}; see Figure~\ref{rootthgs}.
\nomenclature[2a]{$\vrho$}{root corner for a general map}%
\nomenclature[2b]{$\vtau$}{root corner for a unicellular map}%

An oriented corner gives a local orientation to the vertex to which it is incident and to the face that contains it (as it is homeomorphic to a disk by definition). For an oriented corner~$\vc$, we denote by~$\sigma(\vc)$ the subsequent oriented corner around the vertex incident to~$\vc$, in the orientation given by~$\vc$. We also denote by~$\varphi(\vc)$ the oriented corner subsequent to~$\vc$ in the contour of the face containing~$\vc$, in the orientation given by~$\vc$. We will use the following involution we call a \emph{root flip}: from a map~$\m$ rooted at~$\vrho$, we define the root flipped map~$\bar\m$ by rerooting~$\m$ at the oriented corner $\varphi\circ\sigma^{-1}(\vrho)$; see Figure~\ref{rootthgs}. We let $\opp(\vc)$ be the oriented corner corresponding to the same corner as~$\vc$, with the opposite orientation. We use the classical notation $V(\m)$ to denote the vertex set of a map~$\m$ and we denote by~$d_\m$ the graph metric on~$V(\m)$. This means that $d_\m(u,v)$ is the length of a shortest path linking the vertex~$u$ to the vertex~$v$.
\nomenclature[3a]{$\vc$}{oriented corner}%
\nomenclature[3b]{$\sigma(\vc)$}{subsequent oriented corner around the vertex}%
\nomenclature[3c]{$\varphi(\vc)$}{subsequent oriented corner in the contour of the face}%
\nomenclature[2e]{$\bar\m$}{root flipped map of~$\m$}%
\nomenclature[3d]{$\opp(\vc)$}{oriented corner with the opposite orientation}%
\nomenclature[2c]{$V(\m)$}{vertex set of the map~$\m$}%
\nomenclature[2d]{$d_\m$}{graph metric on $V(\m)$}%

An important combinatorial feature a map can have is to be bipartite: a map is \emph{bipartite} if its vertex set can be partitioned into two subsets such that every edge links a vertex from one subset to a vertex from the other subset. Equivalently, a map is bipartite if and only if every cycle made of edges in the map has an even number of edges\footnote{This is a classical fact in graph theory. If a graph is bipartite, the vertices along any cycle alternate between the two subsets, so the cycle must have even length. Conversely, the parity of the distances to a fixed vertex of the graph gives a satisfactory partition: indeed, if the distances to the fixed vertex of two neighboring vertices had the same parity, we would have an odd-length cycle.}.

\subsection{Aim of the paper}

In this work, we present a bijection between the set of pointed bipartite maps of~$\cS$ and pairs consisting of what we call a \emph{well-labeled unicellular mobile} and a parameter $\eps\in\{+,-\}$. A well-labeled unicellular mobile is a one-face map with green or white vertices whose white vertices carry positive integer labels satisfying certain compatibility relations, which need a bit more background to be properly stated (Definition~\ref{wlum}).

If the surface~$\cS$ is orientable, we recover the orientable generalization of the famous Bouttier--Di Francesco--Guitter bijection \cite{bouttier04pml,chapuy07brm} (see also~\cite{chapuy09ori}), and the following basic properties continue to hold in the nonorientable case (Proposition~\ref{faces}). If $(\m,v^\ooo)$ is a pointed bipartite map and $((\fu,\lab),\eps)$ denotes the corresponding pair, then
\nomenclature[4a]{$(\m,v^\ooo)$}{pointed map}%
\nomenclature[4b]{$\eps$}{parameter in $\{+,-\}$}%
\begin{enumerate}[label=(\textit{\roman*})]
	\item $V(\m)\setminus\{v^\ooo\}$ corresponds to the white vertices of~$\fu$ and the label of a white vertex is given by its distance to~$v^\ooo$ in~$\m$;\label{propi}
	\item the faces of~$\m$ correspond to the green vertices of~$\fu\,$: moreover, the degree of a face of~$\m$ is twice the degree of the corresponding green vertex;\label{propii}
	\item the maps~$\m$ and~$\fu$ have the same number of edges;
	\item for each white corner of~$\fu$, exactly one edge of~$\m$ links the incident vertex to a vertex that is closer to~$v^\ooo$.
\end{enumerate}

Property~\ref{propi} is absolutely crucial from a metric point of view, as the labeled unicellular map somehow captures part of the metric information of the map, namely all the distances to the distinguished vertex~$v^\ooo$.

Our construction is based on a rule that gives an orientation to every corner of the map. We introduce what we call \emph{level loops}; these can be thought of as contour lines in topography, where the height of a given vertex is its distance to the distinguished vertex~$v^\ooo$. The orientation of the root gives a canonical orientation to all these level loops. Using these local orientations, we then apply rules similar to those in the orientable case in order to complete the construction.

The inverse construction uses the same mechanism as in the orientable case. One adds first an extra vertex labeled~$0$ inside the unique face of the map. Then, for each corner of the unicellular map incident to a white vertex, one adds an edge linking it to the \textbf{next} corner incident to a white vertex with label strictly smaller (the vertex possibly being the extra label~$0$-vertex). Finally, one removes the initial edges and only keeps the added edges. The whole point in the nonorientable case is to understand how the labels interact with each other and to define the boldfaced word ``next'' of the previous sentence. Our construction will give, also from the orientation of the root, orientations to all the corners incident to white vertices. The word ``next'' will then be understood from these local orientations.

\begin{figure}[ht!]
		\psfrag{2}[][]{$1$}
		\psfrag{3}[][]{$2$}
	\centering\includegraphics[width=6cm]{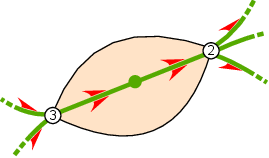}
	\caption{The natural orientation of the corners given by the contour of the unique face cannot work. The half arrowheads indicate the direction of the contour given by the orientation of the root.}
	\label{natoribug}
\end{figure}

As we are working with rooted one-faced maps, the orientation of the root gives a natural orientation to all the corners simply by following the contour of the unique face. Although it would be tempting to perform the mapping with this orientation for the corners, we can see that it would not be satisfactory. Consider for instance a degree~$2$-green vertex whose white neighbors are labeled~$1$ and~$2$ and suppose that, in the contour of the unique face in the vicinity of the green vertex, the vertex labeled~$2$ is visited right before the vertex labeled~$1$ from both sides. We then would expect that both extra edges have to be drawn from the vertex labeled~$2$ to the vertex labeled~$1$. This creates a degree~$2$-face corresponding to the degree~$2$-green vertex, hence violating property~\ref{propii} above; see Figure~\ref{natoribug}. As a result, we will have to deal with two orientations on our maps, these two orientations being the same in the orientable case.

\begin{note}
The orientation process for the level loops we present in this work differs from that of the extended abstract~\cite{bettinelli15fpsac}. It is a bit more complicated to describe but has the benefit of making the encoding maps simpler to describe.
\end{note}

\bigskip

We will then see how to extend our bijection to general maps, which are not necessarily bipartite. In the case of triangulations, that is, maps with only faces of degree~$3$, the encoding mobiles happen to have a particularly simple structure. This allows us to recover the following enumeration result (\cite[Theorem~$1$]{gao91}; see also \cite{eynard12}).

\begin{prop}\label{enum}
The number of (rooted) triangulations of~$\cS$ with $2n$ faces (and thus~$3n$ edges and $n+2-2h$ vertices, by the Euler characteristic formula) is asymptotically equivalent to
$$c_\cS\, n^{5(h-1)/2}\big(12\sqrt 3\big)^n\,,$$
where $h$ is the type of~$\cS$ and $c_\cS$ is a constant that depends on~$\cS$.
\end{prop}
\nomenclature[1c]{$c_\cS$}{asymptotic constant for the number of triangulations of~$\cS$}%

The constants~$c_\cS$ for various~$\cS$ have been the focus of several studies: they are \emph{universal} in the sense that they appear for many different classes of maps~\cite{Gao93pattern} and may be computed through nonlinear recursions~\cite{bender86anr}. Their generating series, properly rescaled, has been shown to satisfy a simple ODE~\cite{GaLeMa08,carrell2014nonorientable}. See also~\cite{Cha19} for an intriguing link with Vorono\"i cell sizes of tessellations of a Brownian surface (in the orientable case). From our bijection, we obtain for~$c_\cS$ a formula involving a summation over cubic maps (Proposition~\ref{cS}). The values for $h\le 1$ are given in Table~\ref{tabcS}.
\begin{table}[ht!]%
\begin{center}
\begin{tabular}{ccc}
\multicolumn{3}{c}{\emph{orientable}}\\[1.5mm]
$h$	& $\cS$	& $c_\cS$\\[.5mm]
\hline\\[-2mm]
$0$	& sphere		& $\dfrac{\sqrt{6}}{\sqrt{\pi}}$\\[4mm]
$1$	& torus		& $\dfrac{1}{8}$
\end{tabular}\qquad\qquad
\begin{tabular}{ccc}
\multicolumn{3}{c}{\emph{nonorientable}}\\[1.5mm]
$h$	& $\cS$	& $c_\cS$\\[.5mm]
\hline\\[-2mm]
$\dfrac12$	& projective plane	& $\dfrac{2^{-3/4}\, 3^{5/4}}{\Gamma(3/4)}$ \\[4mm]
$1$			& Klein bottle		& $\dfrac{3}{2}$
\end{tabular}
\end{center}
\caption{Value of the constant~$c_\cS$ from Proposition~\ref{enum} for $h\le 1$.}
\label{tabcS}
\end{table}

For small values of~$h$, the generating function of triangulations can also be computed: we recover \cite[Theorem~$3$]{gao91} and we add the case of the Klein bottle.
\begin{prop}\label{gfun}
The generating function of triangulations counted with weight~$x$ per vertex is given by
\begin{align*}
&\frac12 \si^3(1-\si)(1-4\si+2\si^2) &&\text{ if~$\cS$ is the sphere,}\\
&\frac12(1-2\si)(1-\si+\si^2)-\frac12\sqrt{1-6\si+6\si^2} &&\text{ if~$\cS$ is the projective plane,}\\
&\frac12\si(1-\si)\big(1-6\si+6\si^2\big)^{-2} &&\text{ if~$\cS$ is the torus,}
\end{align*}
and
$$3\si(1-\si)\big(1-6\si+6\si^2\big)^{-2}\Big(7-30\si+30\si^2-6(1-2\si)\sqrt{1-6\si+6\si^2}\Big)$$
if~$\cS$ is the Klein bottle, where~$\si$ is an algebraic function of~$x$ given by
$$x=\frac12 \si(1-\si)(1-2\si)\qquad \si(0)=0.$$
\end{prop}

See~\eqref{sigmax} for a combinatorial interpretation of~$\si$. In the orientable case, our bijection is the previously known orientable generalization of the Bouttier--Di Francesco--Guitter bijection and Propositions~\ref{enum} and~\ref{gfun} are obtained by the method we use without needing our general bijection. As we did not find this in the literature, even in the easy case of the sphere, we fully treat it here. Note however that, in the case of the sphere, a computation that is similar in spirit was done in~\cite{BDG02census} and, in the orientable case, Chapuy~\cite{chapuy09ori} used the same method in order to obtain similar formulas for maps with face degrees belonging to a given subset of $2\N$ (he also derived results in the more general context of $m$-constellations and $m$-hypermaps). Finally, the analog of Proposition~\ref{enum} was obtained by the same method for quadrangulations on an orientable surface in~\cite{chapuy07brm}, as well as on a nonorientable surface in~\cite{ChDo15bij}. 

Unfortunately, the constraints on the labels of a well-labeled unicellular mobile in general are too intricate to derive similar enumeration results. In fact, the labels satisfy \emph{global} constraints instead of just \emph{local} constraints as soon as the maps we consider are neither triangulations nor quadrangulations. 

\bigskip
The remainder of the paper is organized as follows. We first focus on bipartite maps. We define in Section~\ref{secwlum} the encoding objects we call well-labeled unicellular mobiles and present the inverse bijection. In Section~\ref{secll}, we introduce level loops and present our bijection. Section~\ref{secpfinv} is devoted to the proof that the mappings presented in the previous sections are bijections that are inverse one from another. We then extend our bijection to general maps in Section~\ref{secgen}, focusing on the particular case of triangulations in Section~\ref{sectri}. In Appendix~\ref{secquad}, we come back to the particular case of bipartite quadrangulations by giving an alternate orientation process for the level loops. Finally, in Appendix~\ref{annpol}, we give alternate figures using the useful representation of a unicellular map as a polygon with paired sides.

\begin{ack}
This work is partially supported by Grant ANR-14-CE25-0014 (GRAAL). The author also acknowledges partial support from the Isaac Newton Institute for
Mathematical Sciences where part of this work was conducted. We thank Guillaume Chapuy and Maciej Do{\l}{\k{e}}ga for stimulating discussions and for sharing their paper~\cite{ChDo15bij} with us, as well as Gr\'egory Miermont for encouraging discussions, especially about further developments on the subject of nonorientable Brownian surfaces. We thank \'Elie de Panafieu and Emmanuel Guitter for discussions on the enumeration application. Finally, we thank anonymous referees for interesting comments leading to an improvement of our presentation.
\end{ack}

\section{Well-labeled unicellular mobiles}\label{secwlum}

We now describe in more detail the encoding objects.

\subsection{Encoding objects for quadrangulations}\label{secencquad}

In the case of quadrangulations, an encoding object is particularly simple to describe. It is a so-called \emph{well-labeled unicellular map}, which is a rooted one-face map of~$\cS$, whose vertices carry positive integers that differ by at most~$1$ for neighboring vertices, and with minimal value equal to~$1$.

\subsection{Labeled unicellular mobiles}

\begin{defi}[Labeled unicellular mobile]
A pair~$(\fu,\lab)$ is a \emph{labeled unicellular mobile} if it satisfies the following conditions:
\begin{itemize}
	\item $\fu$ is a rooted one-face map of~$\cS$ whose vertex set is partitioned into $V_\oov(\fu)\sqcup V_\circ(\fu)$ in such a way that every edge links a vertex from $V_\oov(\fu)$ to a vertex from $V_\circ(\fu)$;
	\item $\lab:V_\circ(\fu)\to\N$ is a function with minimum~$1$;
	\item the root vertex belongs to $V_\circ(\fu)$.
\end{itemize}
\end{defi}
\nomenclature[4c]{$(\fu,\lab)$}{labeled unicellular mobile}%
\nomenclature[4d]{$\fu$}{rooted one-face map of~$\cS$ whose vertices are properly colored green or white}%
\nomenclature[4e]{$V_\oov(\fu)$}{set of green vertices of~$\fu$}%
\nomenclature[4f]{$V_\circ(\fu)$}{set of white vertices of~$\fu$ (the one that are labeled)}%
\nomenclature[4g]{$\lab$}{labeling function on $V_\circ(\fu)$}%

The integer~$\lab(v)$ is called the \emph{label} of~$v$; a corner~$c$ or an oriented corner~$\vc$ incident to~$v$ is also said to have \emph{label} $\lab(c)=\lab(\vc)\de\lab(v)$. The elements of~$V_\oov(\fu)$ will be called \emph{green vertices} and those of~$V_\circ(\fu)$ will be called \emph{white vertices}. A corner incident to a white vertex (resp.\ green vertex) is called a \emph{white corner} (resp.\ \emph{green corner}). As~$\fu$ has a unique face, all its corners inherit a canonical orientation from the root~$\vtau$ and are naturally arranged in the order~$\vtau$, $\varphi(\vtau)$, $\varphi^2(\vtau)$, \ldots, $\varphi^{2\#E(\fu)-1}(\vtau)$. This orientation of the corners will be called their \emph{root-induced orientation} and this arrangement will be called the \emph{root-induced contour order}.

\begin{framed}
In what follows, unless explicitly mentioned, we will only consider white corners of~$\fu$, that is, corners that are incident to vertices of~$V_\circ(\fu)$.
\end{framed}

The encoding objects will be labeled unicellular mobile whose labels satisfy extra interaction constraints. In the orientable case, these constraints are the following.

\paragraph{Orientable case.}
If $\cS$ is orientable, the condition of~\cite{bouttier04pml} reads: a labeled unicellular mobile $(\fu,\lab)$ is \emph{well labeled} if, for every corner~$c$, the label of the first subsequent corner is greater than or equal to $\lab(c)-1$. This condition can be restated in the following way, which better fits our definition in the nonorientable case. The labeled unicellular mobile $(\fu,\lab)$ is well labeled if, for every corner~$c$, the label of the first subsequent corner with label strictly smaller than~$\lab(c)$ is equal to $\lab(c)-1$.

The former definition is more compact and simpler to manipulate, especially if one is interested in enumeration, but the latter definition better shows an important property on which the Bouttier--Di~Francesco--Guitter bijection is based. Namely, one obtains from the encoding mobile $(\fu,\lab)$ the edges of the original map by linking each corner of~$\fu$ to the first subsequent corner with strictly smaller label. The vertices of the original map are the vertices of~$V_\circ(\fu)$ (with the addition of one extra vertex) and, in the end, the labeling function gives graph distances to a fixed vertex. As a result, the labels of two neighboring vertices in the original map have to differ by at most one; this may only happen if the labeling is such that the unicellular mobile is well labeled. 

This condition has the advantage of being local, in the sense that one only needs to look at the labels of pairs of subsequent corners in order to see whether a labeled map is well labeled or not. This is fortuitous; in general, it is not to be expected that the interaction constraints become local. This is because they rely on local orientations coming from a global process, which cannot be rendered local, by essence of nonorientability. More precisely, for the phrase ``the first subsequent corner with label strictly smaller than~$\lab(c)$'' to make sense, one needs to orient the corner~$c$. This will be done in accordance with a global orientation process.

\subsection{Coherent orientation of the corners}\label{seccoh}

The idea is that we cannot orient all the corners separately; all the corners of what we call a \emph{corner cycle} should be oriented in accordance. In order to orient the corners of a given corner cycle, we will use the orientation given by the root: more precisely, we will take the first corner of the corner cycle in the root-induced contour order and orient it according to its root-induced orientation, hence forcing the orientation of the whole corner cycle.

To define a corner cycle, we need a starting point, namely a corner. We also need an arbitrary starting orientation: we use the root-induced orientation for convenience but taking the opposite orientation would yield the same result.

Recall the definitions of the mappings~$\sigma$ and~$\varphi$ on oriented corners depicted on Figure~\ref{rootthgs}. Let us now explain how to orient a corner~$c$ of a labeled unicellular mobile $(\fu,\lab)$, constructing its corner cycle in the process. We first temporarily orient~$c$ by letting~$\vc$ be the corner~$c$ oriented according to its root-induced orientation. We then algorithmically construct a list of white and green oriented corners as follows. We initialize the process by letting $\vc_{\textrm{current}}$ be the oriented corner~$\vc$ and by letting the list be empty. Then, iteratively, we add the corner corresponding to~$\vc_{\textrm{current}}$ to the list and,
\begin{itemize}
	\item if $\vc_{\textrm{current}}$ is white, then we update $\vc_{\textrm{current}}$ to the green oriented corner $\varphi(\vc_{\textrm{current}})$;
	\item if $\vc_{\textrm{current}}$ is green, then
	\begin{itemize}
		\item if $\lab\big(\varphi(\vc_{\textrm{current}})\big)\ge\lab(c)$, then we update $\vc_{\textrm{current}}$ to $\varphi(\vc_{\textrm{current}})$;
		\item if $\lab\big(\varphi(\vc_{\textrm{current}})\big)<\lab(c)$, then we update $\vc_{\textrm{current}}$ to~$\sigma(\vc_{\textrm{current}})$. 
	\end{itemize}
\end{itemize}
We stop when $\vc_{\textrm{current}}$ is updated to~$\vc$.

\begin{lem}\label{lemterm}
The above algorithm terminates and outputs a list that does not contain both a corner~$\va$ and its opposite~$\opp(\va)$.
\end{lem}

\begin{pre}
First note that, whenever~$\vc_{\textrm{current}}$ is white, its label is larger than or equal to~$\lab(c)$. Let us consider the first time~$\vc_{\textrm{current}}$ is updated either to an oriented corner that belongs to the current list or to the opposite of such an oriented corner. Let~$\va$ be the value of~$\vc_{\textrm{current}}$ at this time and~$\vb$ be the preceding value of~$\vc_{\textrm{current}}$. As there is a finite number of corners, this time always exists and it suffices to show that $\va=\vc$.

We denote by~$\update$ the function mapping an oriented corner to its update and by~$L$ the list right before the considered time (when~$\vc_{\textrm{current}}=\vb$). We thus have $\va=\update(\vb)$ and it is easy to check that we also have $\opp(\vb)=\update(\opp(\va))$. This implies that $\opp(\va)\notin L$ as this would imply $\opp(\vb)\in L$, a contradiction with the definition of~$\vb$.

As a result, $\va\in L$. Let us argue by contradiction and assume that $\va\neq\vc$. Then there exists $\vb'\neq\vb$ such that $\va=\update(\vb')$. The only possibility is that one of~$\vb$, $\vb'$ is white, say~$\vb'$, and the other is green. Then $\lab\big(\varphi(\vb)\big)<\lab(c)$ and $\va=\varphi(\vb')=\sigma(\vb)$. The latter equality implies that $\lab(\vb')=\lab\big(\varphi(\vb)\big)$, a contradiction with the observation that the label of~$\vc_{\textrm{current}}$ is always larger than or equal to~$\lab(c)$. 
\end{pre}

\begin{figure}[ht!]
		\psfrag{2}[][][.8]{$2$}
		\psfrag{3}[][][.8]{$3$}
		\psfrag{4}[][][.8]{$4$}
	\centering\includegraphics[width=8cm]{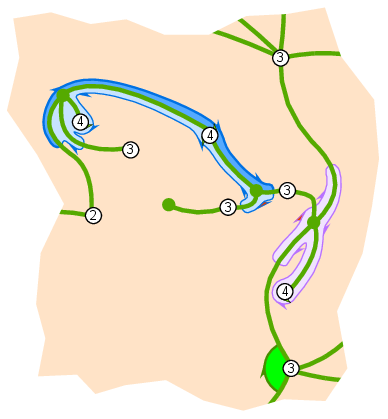}
	\caption{Coherently orienting four corners with label~$4$. There are two corner cycles in the figure. The first corner in the root-induced contour order of each corner cycle is in red. One corner run containing one white corner and two green corners is highlighted at the top of the figure.}
	\label{corcyc}
\end{figure}

Now, consider the smallest integer $i\ge 0$ such that either $\varphi^i(\vtau)$ or $\opp(\varphi^i(\vtau))$ appears in the above output list. If it is~$\varphi^i(\vtau)$ that appears in the list, then we keep for~$c$ its temporary orientation; otherwise, we reverse the orientation of~$c$ and we modify the list by reversing each of its elements, as well as the whole list (by rewriting it backwards); see Figure~\ref{corcyc}. Note that the previous lemma implies that~$\varphi^i(\vtau)$ and~$\opp(\varphi^i(\vtau))$ cannot both belong to the output list.

\begin{rem}
Considering only odd~$i$'s, which amounts to considering only green corners, in the above process produces the same output. This is due to the fact that every white corner that appears in the list is preceded and succeeded by its neighboring green corners.
\end{rem}

\begin{defi}[Coherent orientation]
The above orientation of the corner~$c$ is called its \emph{coherent orientation}.
\end{defi}

\begin{defi}[Corner cycle]\label{defcc} See Figure~\ref{corcyc}.
\begin{itemize}
	\item We call \emph{oriented corner cycle} of~$c$ the above list made into an oriented cycle, that is, considered up to cyclic shift.
	\item The (nonoriented) \emph{corner cycle} of~$c$ is its oriented corner cycle without the data of the orientation, that is, the cycle, considered up to reversing, consisting of the corners corresponding to the oriented corners of the oriented corner cycle.
	\item A \emph{corner run} of a corner cycle is a subset of a corner cycle between two consecutive applications of~$\sigma$ in the above algorithm.
\end{itemize}
\end{defi}

Note that a corner run always starts and ends with a green corner and may be reduced to a single green corner: in this case, it will be called \emph{trivial}. In the algorithm defining the corner cycles, the property that~$\va$ is updated to~$\vb$ if and only if $\opp(\vb)$ is updated to~$\opp(\va)$ yields that the corner cycles do not depend on the orientation of the root, only the \emph{oriented} corner cycles do.

\begin{rem}
The corner cycle of any corner with label~$\lab(c)$ that appears in the corner cycle of~$c$ is the same. As a result, one can actually orient all these corners at once. Furthermore, the corner cycle of a corner labeled~$1$ consists in the contour of the face of~$\fu$: as a result, the corners labeled~$1$ are coherently oriented in their root-induced orientation. 
\end{rem}

\begin{defi}[Successor]\label{defsuc}
Let~$c$ be a corner with label $\lab(c)\ge 2$. Along the face of~$\fu$, in the order given by the coherent orientation of~$c$, the first subsequent corner with label strictly smaller than~$\lab(c)$ is called the \emph{successor} of~$c$; it is denoted by $\suc(c)$.
\nomenclature[6]{$\suc(c)$}{successor of the corner~$c$}%

Formally, if~$\vc$ is the corner~$c$ coherently oriented, then~$\suc(c)$ is the corner corresponding to~$\varphi^i(\vc)$, where
$i\de\min\big\{j\in 2\N : \lab\big(\varphi^j(\vc)\big)<\lab(c)\big\}$.
\end{defi}

We may now define the encoding objects.

\begin{defi}[Well-labeled unicellular mobile]\label{wlum}
A \emph{well-labeled unicellular mobile} is a labeled unicellular mobile such that, for every corner~$c$ with label $\lab(c)\ge 2$, we have $\lab(\suc(c))=\lab(c)-1$.
\end{defi}

\paragraph{Known cases.} Plainly, in the orientable case, Definition~\ref{wlum} generalizes the notion of well-labeled unicellular mobile, as the coherent orientation coincides with the root-induced orientation. In the case of quadrangulations, as explained in the introduction, all the green vertices have degree~$2$ and can thus be discarded. Then Definition~\ref{wlum} is equivalent to the property that the labels of neighboring vertices differ by at most one. Indeed, if this property holds, then the condition of Definition~\ref{wlum} is automatically fulfilled, no matter how the corners are oriented. Conversely, suppose that there are two neighboring vertices whose labels differ by at least~$2$. Then consider an oriented corner~$\vc$ incident to the larger label vertex and such that~$\varphi(\vc)$ is incident to the smaller label vertex. Using the notation~$|\cdot|$ to denote the corner corresponding to an oriented corner, the corner cycle of~$|\vc|$ will thus contain~$|\vc|$ and~$|\sigma(\vc)|$ next to each other. The labels read $\lab(\vc)=\lab(\sigma(\vc))>\lab(\varphi(\vc))+1$, so that either $|\varphi(\vc)|=\suc(|\vc|)$ or $|\si^{-1}\circ\varphi(\vc)|=\suc(|\sigma(\vc)|)$. In both cases, this is in contradiction with Definition~\ref{wlum}.

A similar simplification will also happen in the case of triangulations; see Section~\ref{sectri}.

\subsection{Mapping from well-labeled unicellular mobiles to pointed bipartite maps}\label{secum}

Now that we have the definition of the encoding objects, we can easily describe the inverse of our encoding bijection. It goes as in the orientable case, once every corner is coherently oriented. We consider a well-labeled unicellular mobile $(\fu,\lab)$ and a parameter $\eps\in\{+,-\}$. We add inside the unique face of~$\fu$ a new vertex~$v^\ooo$ with label $\lab(v^\ooo)\de 0$, and we extend the definition of successor by defining the successor of any corner with label~$1$ as the unique corner incident to~$v^\ooo$. Then, for every corner, we add a new edge, called \emph{black edge}, linking it to its successor.

\begin{figure}[ht!]
		\psfrag{l}[][][.8]{$\ell$}
		\psfrag{-}[][][.65]{$\ell{-}1$}
	\centering\includegraphics[width=60mm]{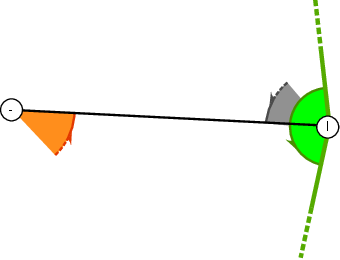}
	\caption{Setting the root of~$\m$, in gray, from the root of~$\fu$, in green. The root of the root flipped map~$\bar\m$ is in orange.}
	\label{rootm}
\end{figure}

We claim that there is a unique way to do so in such a way that the added edges do not intersect with each other or with the edges of~$\fu$, except at their extremities. We consider the embedded graph~$\m$ whose vertex set is $V_\circ(\fu)\cup\{v^\ooo\}$ and whose edges are the added edges. We root it at the oriented corner preceding, in the orientation of the root of~$\fu$, the added edge linking the root corner of~$\fu$ to its successor; see Figure~\ref{rootm}.

\begin{figure}[ht!]
		\psfrag{v}[][]{$v^\ooo$}
		\psfrag{a}[][][.8]{$a$}
		\psfrag{b}[][][.8]{$b$}
		\psfrag{c}[][][.8]{$c$}
		\psfrag{d}[][][.8]{$d$}
		\psfrag{e}[][][.8]{$e$}
		\psfrag{f}[][][.8]{$f$}
		\psfrag{g}[][][.8]{$g$}
		\psfrag{0}[][][.8]{$0$}
		\psfrag{1}[][][.8]{$1$}
		\psfrag{2}[][][.8]{$2$}
		\psfrag{3}[][][.8]{$3$}
		\psfrag{4}[][][.8]{$4$}
	\centering\includegraphics[width=.95\linewidth]{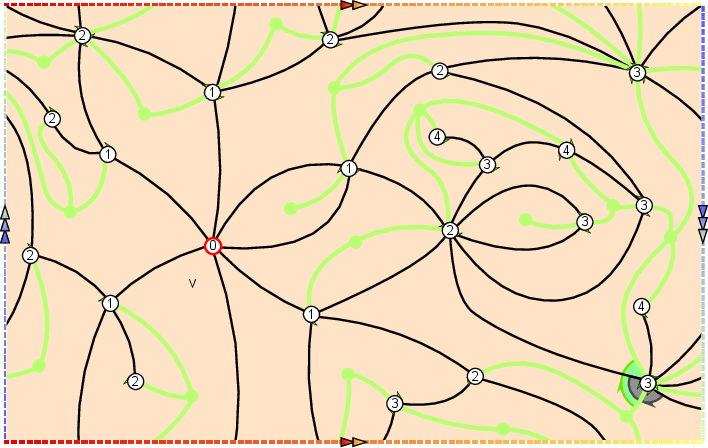}
	\caption{The bijection, from a well-labeled unicellular mobile to a pointed bipartite map (on the Klein bottle). See also Figure~\ref{bij_um_p} for a polygon representation of the same figure where the oriented corner cycles are depicted.}
	\label{bij_um}
\end{figure}

We set the output of our mapping to be $\Psi\big((\fu,\lab),-\big)\de(\m,v^\ooo)$ and $\Psi\big((\fu,\lab),+\big)\de(\bar\m,v^\ooo)$, where~$\bar\m$ denotes the root flipped version of~$\m$ (we will see in Proposition~\ref{propmu} that~$\m$ is actually a map). Figures~\ref{bij_um} and~\ref{bij_um_p} show an example of the mapping.
\nomenclature[9d]{$\Psi$}{mapping $\cU\times\{+,-\} \to \cB^\ooo$}

\bigskip
It remains to prove the claim that the new edges may be drawn in a noncrossing fashion. Let us consider a corner~$c$ with label $\lab(c)\ge 2$ and the chain of edges of~$\fu$ in the contour of the face of~$\fu$, from~$c$ to~$\suc(c)$ in the coherent orientation of~$c$. Then the edge linking~$c$ to~$\suc(c)$ must be drawn in such a way that~$v^\ooo$ lies outside the disk delimited by the above chain and the added edge; this is because there is at least one corner labeled~$1$ outside of this region. Now, with this extra constraint, there is a unique way to add the desired edges as long as the corners to be linked lie in the same face of the map consisting of~$\fu$ together with the previously added edges. 

Let us then argue by contradiction and assume that the corner~$c$ is to be linked to the corner~$\suc(c)$ (possibly incident to~$v^\ooo$) and that they do not lie in the same face. If $\lab(c)\ge 2$, we consider the interval~$I$ of corners of~$\fu$ from~$c$ to~$\suc(c)$ in the coherent orientation of~$c$ ($c$ and~$\suc(c)$ excluded). Then there is a previously added edge linking a corner~$d$ of~$I$ to a corner~$d'$ outside of~$I\cup\{c,\suc(c)\}$. As all the corners of~$I$ have label larger than~$\lab(c)$, we have $\lab(d)\ge\lab(c)$, and, by construction, $\lab(d')=\lab(d)\pm 1$. 
\begin{itemize}
	\item If $\lab(d')=\lab(d)+ 1$, then~$d$ is, in the coherent orientation of~$d'$, the first subsequent corner of~$d'$ with label strictly smaller than~$\lab(d')$. This is not possible as either~$c$ or~$\suc(c)$ arrive before~$d$ in this order and they both satisfy the label requirement.
	\item If $\lab(d')=\lab(d)- 1$, then~$d'=\suc(d)$. This can only happen if $\lab(d)=\lab(c)$ and~$d$ is coherently oriented in the opposite direction as the coherent orientation of~$c$. But this is not possible, as~$d$ would thus belong to the corner cycle of~$c$.
\end{itemize}
Finally, if $\lab(c)=1$, then there is a previously added edge linking a corner~$d$ to a corner~$d'$ such that, in the coherent orientation of~$d$, the first subsequent corner of~$d$ with label strictly smaller than~$\lab(d)$ is~$d'$, and~$c$ belongs to the interval of corners between~$d$ and~$d'$. This is impossible as~$c$ has minimum label. This proves the claim.

\begin{prop}\label{propmu}
The output of the above mapping $\Psi\big((\fu,\lab),\pm\big)$ is a pointed bipartite map.
\end{prop}

\begin{pre}
It will be enough to show that~$\m$ is a bipartite map, as this immediately implies that $\bar\m$ is also a bipartite map.

Let us consider two adjacent corners~$c$ and~$c'$ of~$\fu$, that is, two subsequent white corners in the contour of the face of~$\fu$, assuming $\lab(c)\ge \lab(c')$. Let also~$\vc$ be the oriented corner corresponding to~$c$ oriented toward~$c'$. Then, the first $\lab(c)-\lab(c')-1$ iterate successors of~$c$ are necessarily coherently oriented in the opposite direction as~$\vc$ and the coherent orientation of the $(\lab(c)-\lab(c'))$-th iterate successor~$c''$ of~$c$ is the same as the coherent orientation of~$c'$. The first observation comes from the fact that the labels of these successors are strictly larger than~$\lab(c')$ and the second observation is due to the fact that~$c''$ and~$c'$ belong to the same corner cycle. As a result, either~$c''=c'$ or $c''\neq c'$ and $\suc(c'')=\suc(c')$ (if $\lab(c)=\lab(c')$, we are necessarily in the second case). Considering the black edges linking the corners to their successors, we thus obtain a chain of black edges linking any two adjacent corners.

We now consider a connected component of the complement of~$\m$. It consists of a finite union of connected components delimited by two green edges incident to the same green vertex and by the chain of black edges considered above. These components are glued together at the green edges incident to the same green vertex and we obtain that the considered connected component is homeomorphic to a disk and contains only one green vertex. As a result, $\m$ is a map and it is bipartite as every black edge links a vertex with odd label to a vertex with even label.
\end{pre}

\section{Level loops}\label{secll}

Throughout this section, we fix a pointed \textbf{bipartite} map $(\m,v^\ooo)$.

\subsection{Definition}\label{secdefll}

The notion of \emph{level loop} plays a crucial part in this work. The level loops of $(\m,v^\ooo)$ lie on the sides of its edges and may cross some edges near their extremities. First, we define the labeling function $\lab:V(\m)\to\Z_+$ by
\begin{equation}\label{defl}
\lab(v)\de d_\m(v^\ooo,v), \qquad v\in V(\m)\,,
\end{equation}
and we call \emph{label} of~$v$ the integer~$\lab(v)$. We extend the definition to a corner~$c$ and to an oriented corner~$\vc$ by setting their label to be the label of the incident vertex, denoted by~$\lab(c)$ and~$\lab(\vc)$. Note that, as the map~$\m$ is bipartite by hypothesis, the labels of two neighboring vertices differ by exactly~$1$. It will be convenient to canonically orient the edges of~$\m$ in such a way that the label of the head is one less than the label of the tail: this orientation will be called the \emph{geodesic orientation}.

We pick an arbitrary oriented corner of~$\m$ and denote by~$i$ its label. We start from this oriented corner and look at the incident edge:
\begin{itemize}
	\item if it leads to a vertex with label smaller than or equal to~$i$, we move along its side up to the subsequent corner of the face;
	\item if it leads to a vertex labeled $i+1$, we cross it.
\end{itemize}
We then iterate the process (with the same~$i$), in the sense that, whenever we encounter an edge, we either move along its side if it leads to a vertex of label smaller than or equal to~$i$ or we cross it if it leads to a vertex labeled $i+1$, until we come back to the initial corner and close the loop.

Plainly, the resulting loop does not depend on the orientation of the initial corner.

\begin{defi}[Level loop]
We call \emph{level loop} issued from the corner~$c$ the loop defined by the above process with~$c$ as the initial corner. Its \emph{level} is the nonnegative integer~$\lab(c)$.
\end{defi}

\begin{figure}[ht!]
		\psfrag{v}[][]{$v^\ooo$}
		\psfrag{0}[][][.8]{$0$}
		\psfrag{1}[][][.8]{$1$}
		\psfrag{2}[][][.8]{$2$}
		\psfrag{3}[][][.8]{$3$}
		\psfrag{4}[][][.8]{$4$}
	\centering\includegraphics[width=.95\linewidth]{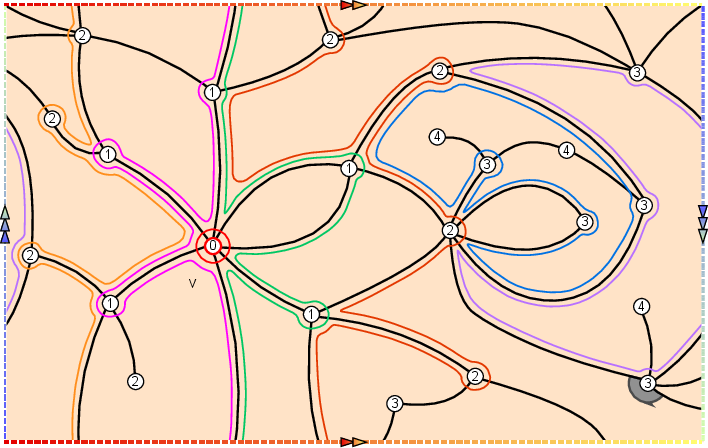}
	\caption{Level loops of a pointed bipartite map. For each face, we did not represent the maximal level loop, which follows its contour. The red loop is at level~$0$, the pink and green loops are at level~$1$, the dark red and orange loops are at level~$2$ and the blue and purple loops are at level~$3$.}
	\label{ll}
\end{figure}

Let us state some immediate properties about level loops. First, the same level loop cannot visit the same corner twice, and two different level loops at same level cannot both visit the same corner. Next, the labels of the corners visited by a loop at level~$i$ are all smaller than or equal to~$i$, and at least one is equal to~$i$. There is only one loop that circles around a vertex: it is the level loop at level~$0$ and it circles around~$v^\ooo$. Indeed, if a loop circles around a vertex labeled~$i$, it means that the vertex is only connected to vertices labeled $i+1$; the only such vertex is~$v^\ooo$. A level loop that does not cross any edge follows the contour of a face; its level is the maximal label among those of the corners of the face. Such loops will be called \emph{maximal level loops}. Figure~\ref{ll} shows all the level loops except the maximal level ones of a pointed bipartite map.

\begin{prop}\label{lltrivial}
Let~$c$ be a corner of a face~$f$. Let~$i$ denote the label of~$c$ and $j\ge i$ denote the maximum label of the corners of~$f$. Then the corner~$c$ is visited by exactly $j-i+1$ different level loops, one for each level from $\{i,i+1,\ldots,j\}$.
\end{prop}

\begin{note}
In the figures, the level loops with higher level will always be drawn further away from the edges than the loops with lower level. This ensures that they do not intersect, since, whenever two loops encounter the same edge, if the lower level one moves along its side, then so does the higher level one. 
\end{note}

\begin{figure}[ht!]
		\psfrag{2}[][][.8]{$2$}
		\psfrag{3}[][][.8]{$3$}
		\psfrag{4}[][][.8]{$4$}
		\psfrag{5}[][][.8]{$5$}
		\psfrag{6}[][][.8]{$6$}
		\psfrag{l}[][][.8]{$\ell$}
		\psfrag{+}[][][.45]{$\ell{+}1$}
		\psfrag{-}[][][.45]{$\ell{-}1$}	
	\centering\includegraphics[width=11cm]{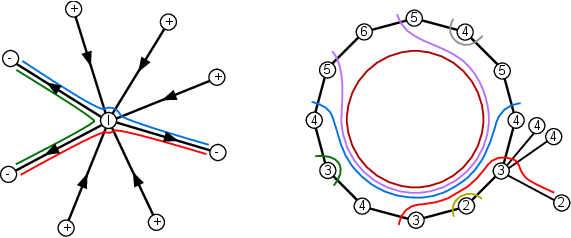}
	\caption{\textbf{Left.} Parts of level loops at level~$\ell$ around a vertex with label~$\ell$. The geodesic orientation of the edges is represented by full arrowheads. \textbf{Right.} Parts of level loops in a face. For instance, the blue loop is at level~$4$; it crosses edges whose extremities are labeled~$4$ and~$5$, and it moves along edges with extremities both labeled~$4$ or less. Note that two parts of level loops at the same level may meet outside the face and thus belong to the same loop.}
	\label{lltyp}
\end{figure}

Figure~\ref{lltyp} shows the situation around a typical vertex and inside a face.

\begin{pre}[Proof of Proposition~\ref{lltrivial}]
Let us consider a level loop at level~$k$ that visits the corner~$c$. As a level loop only visits corners whose labels are smaller than its level, we have $k\ge i$. Moreover, if $k>j$, then the definition forces the loop to follow the contour of~$f$. This is not possible as this loop would not visit any corner labeled~$k$, thus contradicting the definition. As a result, there is at most one level loop that visits~$c$ for each level from $\{i,i+1,\ldots,j\}$.

Now, let $k\in\{i,i+1,\ldots,j\}$. Let us arbitrarily orient~$c$ and consider the next corner with label~$k$ in the contour of~$f$. Then the level loop issued from this corner is at level~$k$ and visits~$c$.
\end{pre}

\begin{figure}[ht!]
		\psfrag{r}[B][B]{$\vrho$}		
		\psfrag{v}[][]{$v^\ooo$}
		\psfrag{0}[][][.8]{$0$}
		\psfrag{1}[][][.8]{$1$}
		\psfrag{2}[][][.8]{$2$}
		\psfrag{3}[][][.8]{$3$}
		\psfrag{4}[][][.8]{$4$}
		\psfrag{a}[Br][Br][.8]{$\pr(\vrho)$}
		\psfrag{b}[Bl][Bl][.8]{$\pr^2(\vrho)$}
		\psfrag{c}[Bl][Bl][.8]{$\pr^3(\vrho)$}
	\centering\includegraphics[width=10cm]{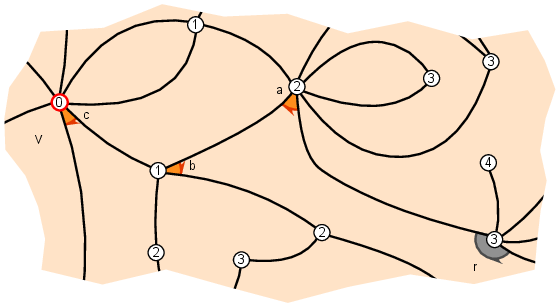}
	\caption{Iterative parent corners of the oriented corner~$\vrho$.}
	\label{parent}
\end{figure}

Let us consider an oriented corner~$\vc$ not incident to~$v^\ooo$. Let us also consider the level loop issued from~$\vc$, oriented in the orientation given by~$\vc$. The first oriented corner with label $\lab(\vc)-1$ after~$\vc$ visited by the loop will be of interest: we call it the \emph{parent corner} of~$\vc$ and denote it by~$\pr(\vc)$; see Figure~\ref{parent}.
\nomenclature[8]{$\pr(\vc)$}{parent corner of~$\vc$}%

\subsection{Orienting the level loops}\label{secori}

We need a canonical way to orient all the level loops. Through the bijection, the level loops will correspond to the corner cycles of the encoding unicellular mobile. Our orientation process for the level loops thus has to match the orientation process of the corner cycles of the encoding mobile. There are several options; the one we present in this work corresponds to the coherent orientation of the corners described in Section~\ref{seccoh}. We proposed in~\cite{bettinelli15fpsac} a different option: it was simpler to describe but its counterpart on the encoding objects was harder to comprehend. We will present the latter option in the particular case of bipartite quadrangulations in Appendix~\ref{secquad}. Before explaining in details the orientation process, let us give some insights that should help the reader.

In regard to the way we defined the coherent orientation, we will need to explore the unique face of the encoding map in its root-induced contour order. We will actually build this cyclic order on the fly while constructing the encoding unicellular map. We will orient the loops one after the other. The orientation of a loop will allow us to ``break'' the loop into several pieces by adding on it marks called \emph{stops}, which will correspond to the edges of the encoding mobile. We will then need to go from a stop to the next one in what corresponds to the root-induced contour order of the encoding mobile. A portion of level loop in between consecutive stops will correspond to a run of successive corners that appears in the corner cycle corresponding to the level loop. The way this run of corners is to be explored is dictated by the global orientation of the face of the encoding mobile and not by the orientation of the level loop.

\bigskip
We break the orientation process into two main operations.

\begin{figure}[ht!]
		\psfrag{1}[][][.8]{$1$}
		\psfrag{2}[][][.8]{$2$}
		\psfrag{3}[][][.8]{$3$}
		\psfrag{4}[][][.8]{$4$}
	\centering\includegraphics[width=6cm]{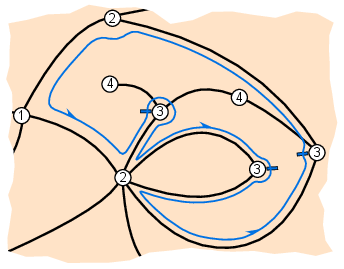}
	\caption{Definition of the stops (symbolized by small rectangles of the same color as the loop) of a level loop from its orientation. This loop is at level~$3$; it thus has stops at the corners labeled~$3$ directly preceded by corners labeled~$2$.}
	\label{stops}
\end{figure}

\paragraph{Definition of stops after orientation [Figure~\ref{stops}].} Whenever we orient a level loop at level~$i$, we add on it stops at the locations of the corners labeled~$i$ the loop visits right after visiting a corner labeled $i-1$. Equivalently, for each face visited by the loop, there is a stop at each corner labeled~$i$ except for the first one.

\begin{figure}[ht!]
		\psfrag{1}[][][.8]{\textit{(a)}}
		\psfrag{2}[][][.8]{\textit{(b)}}
		\psfrag{3}[][][.8]{\textit{(c)}}
		\psfrag{4}[][][.8]{\textit{(d)}}
		
		\psfrag{l}[][][.8]{}
		\psfrag{+}[][][.5]{}
		\psfrag{-}[][][.5]{}
		
		\psfrag{d}[Bl][Bl][.7]{$\si^j(\vc)$}
		\psfrag{e}[Br][Br][.7]{$\varphi^{-k}\circ\si^j(\vc)$}
		\psfrag{f}[Bl][Bl][.7]{$\upd(\vc)$}
		\psfrag{a}[B][B][.7]{\textcolor{red}{$\lambda$}}
		\psfrag{c}[][][.7]{$\vc$}
	\centering\includegraphics[width=.95\linewidth]{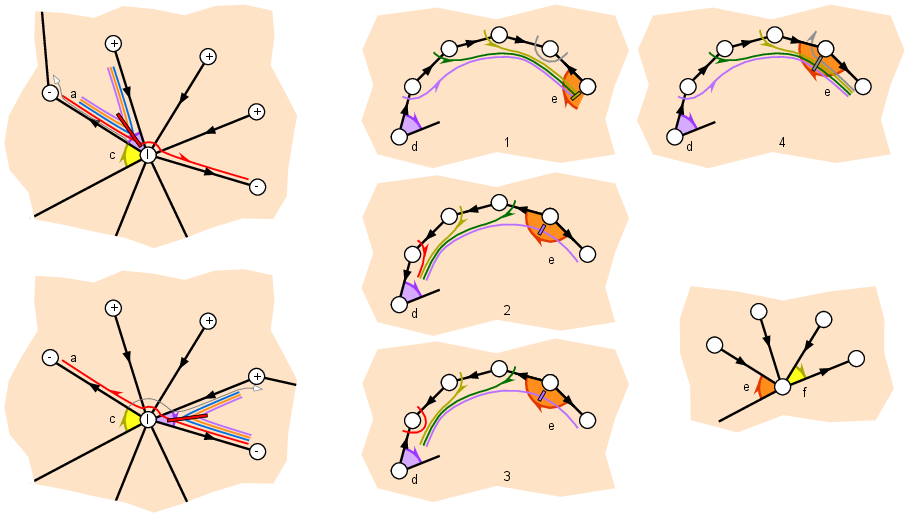}
	\caption{Updating the current corner. All the possibilities are represented. \textbf{Left.} We turn around the initial vertex until we reach a stop (the light gray arrow indicates where the process continues). On the top, the loop~$\lambda$ may or may not have been oriented before. On the bottom, it has been oriented before. \textbf{Middle and top right.} We follow the contour of the face until we reach a stop. The orientations of the gray loop in~\textit{(a)} and the red loop in~\textit{(c)} do not intervene at this stage; they may or may not have been oriented before. In~\textit{(b)} and~\textit{(c)}, the purple loop may be oriented both ways, either at this stage or before (if at this stage, necessarily in the orientation of~$\sigma^j(\vc)$). In~\textit{(d)}, the gray loop may have been oriented before or at this stage. \textbf{Bottom right.} We turn around the final vertex until we reach an edge that brings us closer to~$v^\ooo$.}
	\label{update}
\end{figure}

\paragraph{Update of the current oriented corner [Figure~\ref{update}].} This step will allow to go from a white corner of the encoding unicellular mobile to the subsequent one in the root-induced contour order, orienting several level loops along the way. We break it into three substeps. We consider an oriented corner~$\vc$ such that $\lab(\varphi(\vc))<\lab(\vc)$.
\begin{description}[align=left, itemindent=-5mm]
	\item[Turning around the initial vertex.] If the level loop~$\lambda$ issued from~$\si(\vc)$ has not already been oriented, then we orient it in the orientation of~$\si(\vc)$ and define its stops as explained above. Then we consider the smallest \textbf{positive} integer~$j$ such that~$\lambda$ has a stop at the location of~$\si^j(\vc)$. (Note that $j=1$ if the level loop was not oriented before.) Next, we orient, in the orientation of~$\si^j(\vc)$, each level loop visiting~$\si^j(\vc)$ that has not already been oriented, and we define its stops. 
	\item[Moving along the contour of the face.] We consider the smallest \textbf{positive} integer~$k$ such that 
	\begin{itemize}
		\item either there is a stop at the location of $\varphi^{-k}\circ\si^j(\vc)$;
		\item or the level loop issued from $\varphi^{-k}\circ\si^j(\vc)$ has not yet been oriented and orienting it in the orientation of $\varphi^{-k}\circ\si^j(\vc)$ yields a stop at the location of~$\varphi^{-k}\circ\si^j(\vc)$.
	\end{itemize}
	In the latter case (which may only happen in case~\textit{(d)} of Figure~\ref{update}), we orient the loop issued from $\varphi^{-k}\circ\si^j(\vc)$ in the orientation of $\varphi^{-k}\circ\si^j(\vc)$ and define its stops.
	\item[Turning around the final vertex.] Finally, we consider the smallest \textbf{nonnegative} integer~$l$ such that $\lab\big(\varphi\big(\si^{l}\circ\varphi^{-k}\circ\si^j(\vc)\big)\big)=\lab\big(\si^{l}\circ\varphi^{-k}\circ\si^j(\vc)\big)-1$.
\end{description}
 We define
\[
\upd(\vc)\de\si^{l}\circ\varphi^{-k}\circ\si^j(\vc)\,.
\]
\nomenclature[8b]{$\upd(\vc)$}{oriented corner into which~$\vc$ is updated in the process orienting the level loops}%

We may now define the orientation process. We initialize it as follows.

\begin{figure}[ht!]
		\psfrag{r}[B][B]{$\vrho$}		
		\psfrag{v}[][]{$v^\ooo$}
		\psfrag{0}[][][.8]{$0$}
		\psfrag{1}[][][.8]{$1$}
		\psfrag{2}[][][.8]{$2$}
		\psfrag{3}[][][.8]{$3$}
		\psfrag{4}[][][.8]{$4$}
		\psfrag{a}[Br][Br][.8]{}
		\psfrag{b}[Bl][Bl][.8]{}
		\psfrag{c}[Bl][Bl][.8]{}
	\centering\includegraphics[width=10cm]{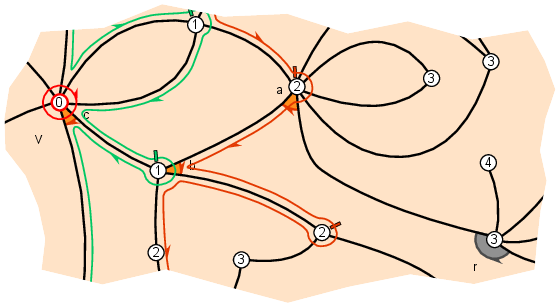}
	\caption{Orienting the level loops issued from the iterate parents of the root.}
	\label{oriparents}
\end{figure}

\paragraph{Initialization [Figure~\ref{oriparents}].} Before we start, we look at the labels of the extremities of the root edge: if the root vertex has a smaller label than the other extremity of the root edge, then we replace~$\m$ by its root flipped version~$\bar\m$.

Then, for each $i\in\{1,\ldots,\lab(\vrho)\}$, we orient the level loop issued from the $i$-th iterate parent~$\pr^i(\vrho)$ of the root in the orientation of~$\pr^i(\vrho)$.

\begin{figure}[ht!]
		\psfrag{v}[][]{$v^\ooo$}
		\psfrag{0}[][][.8]{$0$}
		\psfrag{1}[][][.8]{$1$}
		\psfrag{2}[][][.8]{$2$}
		\psfrag{3}[][][.8]{$3$}
		\psfrag{4}[][][.8]{$4$}
		\psfrag{a}[][][.8]{$a$}
		\psfrag{b}[][][.8]{$b$}
		\psfrag{c}[][][.8]{$c$}
		\psfrag{d}[][][.8]{$d$}
		\psfrag{e}[][][.8]{$e$}
		\psfrag{f}[][][.8]{$f$}
		\psfrag{g}[][][.8]{$g$}
	\centering\includegraphics[width=.95\linewidth]{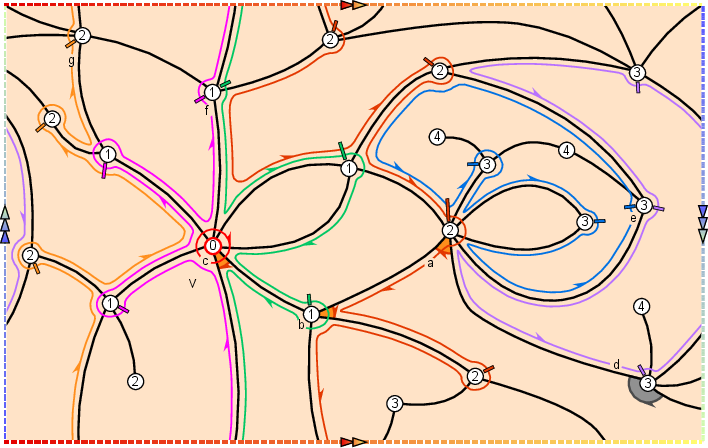}
	\caption{Orientation of the level loops. Recall that the maximal level loops are not represented (this bears no effects). The loops are oriented one after the other, in the order given by the letters a, b, c, d, e, f, g.}
	\label{orill}
\end{figure}

\paragraph{Orientation process [Figure~\ref{orill}].} 
After performing the initialization step, we let the current oriented corner be the root and we iteratively update it using the update step, until every level loop is oriented.

\begin{prop}\label{oriterm}
The above orientation process terminates: it assigns an orientation to all the level loops.
\end{prop}

The proof is postponed to Section~\ref{secpfinv}.

\subsection{Mapping from pointed bipartite maps to well-labeled unicellular mobiles}\label{secmu}

We may now present our bijection: we apply to the pointed bipartite map $(\m,v^\ooo)$ the following steps.

\paragraph{Step 1. Labeling the vertices and setting the value of the parameter~$\eps$.}\refstepcounter{step}\label{step1}
We define the labeling function~$\lab$ by~\eqref{defl}. 
If the root vertex of~$\m$ has a larger label than the other extremity of the root edge, we set $\eps\de -\,$; otherwise, we replace~$\m$ by its root flipped version~$\bar\m$ and set $\eps\de +$.

\paragraph{Step 2. Constructing and orienting the level loops.}\refstepcounter{step}\label{step2}
We construct the level loops and orient them by the process of Section~\ref{secori}, keeping track of the stop locations.

\paragraph{Step 3. Linking the stop corners to green vertices.}\refstepcounter{step}\label{step3}
We add an extra vertex in the middle of each face of~$\m$: these vertices will be called \emph{green vertices}, in contrast with the original vertices of~$\m$, which we will call \emph{white vertices}. Inside each face, we link in a noncrossing fashion by \emph{green edges} the green vertex of the face to all the corners of the face where there is a stop.

\begin{figure}[ht!]
		\psfrag{l}[][][.8]{$\ell$}
		\psfrag{-}[][][.65]{$\ell{-}1$}
	\centering\includegraphics[width=40mm]{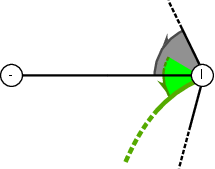}
	\caption{Rooting~$\fu$ from the root of~$\m$. The root of~$\m$ is represented in gray and the root of~$\fu$ is in green. The existence of the green edge (of~$\fu$) is ensured by the rules of the construction.}
	\label{rootu}
\end{figure}

\begin{figure}[ht!]
		\psfrag{v}[][]{$v^\ooo$}
		\psfrag{0}[][][.8]{$0$}
		\psfrag{1}[][][.8]{$1$}
		\psfrag{2}[][][.8]{$2$}
		\psfrag{3}[][][.8]{$3$}
		\psfrag{4}[][][.8]{$4$}
	\centering\includegraphics[width=.95\linewidth]{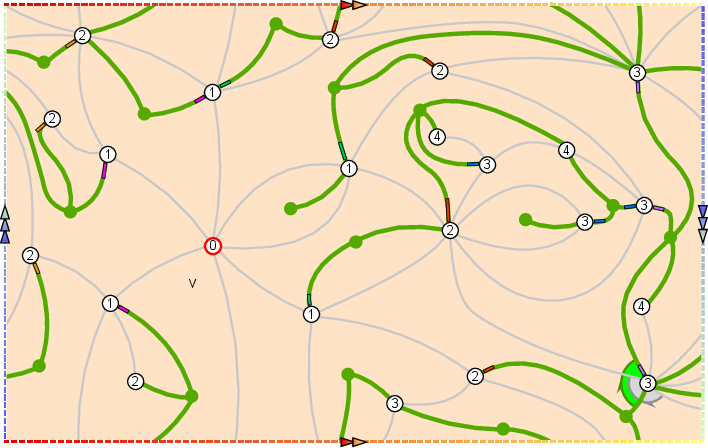}
	\caption{The bijection, from a pointed bipartite map to a well-labeled unicellular mobile. The edges of the original map have been grayed out. The stops of the not represented maximal level loops are at all the corners with maximal label among their face.}
	\label{bij_mu}
\end{figure}

\paragraph{Step 4. Discarding original edges and rooting the resulting map.}\refstepcounter{step}\label{step4}
Finally, we consider the embedded graph~$\fu$ whose vertex set consists of the union of $V_\circ(\fu)\de V(\m)\setminus\{v^\ooo\}$ with the set $V_\oov(\fu)$ of green vertices and whose edge set is composed of the green edges. (Note that this definition makes sense: there are no green edges with extremity~$v^\ooo$, as there are no stops incident to~$v^\ooo$.) We root it with the convention depicted on Figure~\ref{rootu}. The white vertices of~$\fu$ inherit the labels from the function~$\lab$: we set $\Phi(\m,v^\ooo)\de ((\fu,\lab),\eps)$; see Figure~\ref{bij_mu}.
\nomenclature[9c]{$\Phi$}{mapping $\cB^\ooo \to \cU\times\{+,-\}$}

Note that, as with the inverse mapping of Section~\ref{secum}, it is not clear at this stage that the output is even a map; this will be shown in Section~\ref{secpfinv}.

\section{The previous mappings are inverse one from another}\label{secpfinv}

For the surface~$\cS$ we consider, we denote by~$\cB^\ooo$ the set of pointed bipartite maps and by~$\cU$ the set of well-labeled unicellular mobiles (Definition~\ref{wlum}). Recall that we denoted by~$\Phi$ the mapping from Section~\ref{secmu} and by~$\Psi$ the one from Section~\ref{secum}. We have not yet established that~$\Phi$ is well defined, as Proposition~\ref{oriterm} has not yet been proved. However, from a pointed bipartite map $(\m,v^\ooo)\in\cB^\ooo$, we may always define a labeled embedded graph~$(\fu,\lab)$ by the construction of Section~\ref{secmu}, replacing the orientation process of Step~\ref{step2} by an arbitrary orientation process.
\nomenclature[9a]{$\cB^\ooo$}{set of pointed bipartite maps of~$\cS$}%
\nomenclature[9b]{$\cU$}{set of well-labeled unicellular mobiles of~$\cS$}%

\begin{prop}\label{faces}
Let $(\m,v^\ooo)\in\cB^\ooo$ and $(\fu,\lab)$ be defined as above from an arbitrary orientation process. Then
\begin{enumerate}[label=(\textit{\roman*})]
	\item $V(\m)=V_\circ(\fu)\sqcup \{v^\ooo\}$ and, for $v\in V_\circ(\fu)$, $\lab(v)=d_\m(v,v^\ooo)\,$;\label{faces_i}
	\item the faces of~$\m$ correspond to $V_\oov(\fu)\,$: moreover, the degree of a face of~$\m$ is twice the degree of the corresponding vertex in $V_\oov(\fu)$;\label{faces_ii}
	\item the map~$\m$ and the embedded graph~$\fu$ have the same number of edges;\label{faces_iii}
	\item with the edges of~$\m$ oriented according to the geodesic orientation, among the edges of~$\m$ in a white corner of~$\fu$, exactly one is outgoing.\label{faces_iv}
\end{enumerate}
\end{prop}

\begin{pre}
Statement~\ref{faces_i} and the fact that the faces of~$\m$ correspond to the green vertices of~$\fu$ are direct consequences of the construction. Let us consider a face of~$\m$ and denote by~$2p$ its degree. We see the face as a $2p$-gon and pair its sides as follows. Let us pick a side and denote by~$i$ and $i-1$ the labels of its extremities. We travel along the boundary from the vertex labeled~$i$ to the vertex labeled $i-1$, and we keep traveling until we successively encounter a corner labeled $i-1$ and a corner labeled~$i$. We match the side we picked with the side linking these corners (note that these two sides are necessarily distinct and that the matching does not depend on which of the two sides we first picked). This gives a perfect matching of the sides of our polygon.

We then consider a loop at level~$i$ that visits the face. Recall that this loop has stops at the corners with label~$i$ that are immediately preceded by a corner labeled~$i-1$. Then, by construction of the loop, if it has a stop at a corner labeled~$i$ of our face, then the previous corner labeled $i-1$ that the loop visits also belongs to our face; we associate with the stop the side linking these two corners. Now observe that two paired sides with labels \edg{$i$-$1$}{$i$} are visited by the same level loop at level~$i$, one from label~$i$ to label $i-1$ and the other one from label $i-1$ to label~$i$, no matter what the orientation of the loop is. As a result, exactly one of these two sides is associated with a stop. As there are~$p$ pairs of sides, there are also~$p$ stops and the green vertex of the face has degree~$p$.

Point~\ref{faces_iii} is then a direct consequence of~\ref{faces_ii}: the number of edges of~$\m$ is half the sum of the degrees of its faces and the number of edges of~$\fu$ is equal to the sum of the degrees of its green vertices. It is also a straightforward consequence of~\ref{faces_iv} below.

Let us show~\ref{faces_iv}. If we consider a vertex of~$\m$ and its neighbors as on the left of Figure~\ref{lltyp}, we see that there will be exactly one stop between two consecutive outgoing edges (the location depending on the orientation of the corresponding loop). As a white corner of~$\fu$ corresponds to the area between two consecutive stops around a vertex of~$\m$, we see that there is always exactly one outgoing edge in such an area.
\end{pre}

\begin{lem}\label{lemul}
Let $(\m,v^\ooo)\in\cB^\ooo$ and $(\fu,\lab)$ be defined as above from an arbitrary orientation process. Then $(\fu,\lab)$ is a labeled unicellular mobile.
\end{lem}

\begin{pre}
The only thing to show is that the embedded graph~$\fu$ is a unicellular map. Adapting arguments from~\cite{chapuy07brm}, we consider the map~$\fM$ whose vertex set is $V(\m)\cup V_\oov(\fu)$ and whose edges are the edges of~$\m$ together with the edges of~$\fu$. We add inside each face of this map a blue vertex and, for each edge of~$\m$, we add a blue dual edge linking the blue vertices of the two incident faces\footnote{The embedded blue graph we obtain is a generalization of the \emph{dual exploration graph} appearing in~\cite{ChDo15bij} (and also represented in blue in the figures of this reference).}. Let us see that this embedded blue graph has only one cycle, which turns around~$v^\ooo$. To this end, let us label each blue edge by the minimal label of the extremities of the black edge it crosses. Let us consider a blue cycle and denote by~$m$ the maximal label of its edges. Let us consider a blue edge~$e$ of the cycle with label~$m$ and a neighboring edge~$e'$ of the cycle; see Figure~\ref{unicel}.

Let~$\tilde e$ and~$\tilde e'$ denote the black dual edges of~$e$ and~$e'$. These edges are incident to the same face of~$\m$, which is split into two connected components by $e\cup e'$; we denote by~$v$ the extremity of~$\tilde e$ labeled~$m$ and by~$v'$ the extremity of~$\tilde e'$ incident to the same connected component as~$v$. As all the corners whose label is a local maximum in the contour of a face are stops of the loop issued from them, the green vertex must belong to the connected component that does not contain~$v$ (recall that the label of~$e'$ is smaller than or equal to~$m$). As a result, there are no stops in the connected component containing~$v$ and the labels of the corners in this component are thus not local maximums in the contour of the face. In particular, if $v'\neq v$, as the label of~$e'$ is smaller than or equal to~$m$, the neighboring vertex of~$v$ must be labeled $m-1$. We conduct the same reasoning with the blue edge~$e''$ of the cycle that intersects the other extremity of~$e$ and denote by~$\tilde e''$ and~$v''$ the counterparts of~$\tilde e'$ and~$v'$. We see that, if~$v$, $v'$ and~$v''$ are three pairwise distinct vertices, there must exist a level loop at level~$m$ that crosses~$\tilde e$. This loop should have a stop at a corner labeled~$m$ in one of the two faces, which is impossible. As a result, either $v'=v$ or $v''=v$.

\begin{figure}[ht!]
		\psfrag{e}[][][.8]{\textcolor{blue}{$e$}}
		\psfrag{t}[][][.8]{$\tilde e$}
		\psfrag{f}[][][.8]{\textcolor{blue}{$e'$}}
		\psfrag{g}[][][.8]{$\tilde e'$}
		\psfrag{h}[][][.8]{\textcolor{blue}{$e''$}}
		\psfrag{v}[][][.8]{$v$}
		\psfrag{w}[][][.8]{$v'$}
		\psfrag{x}[][][.8]{$v''$}
		\psfrag{m}[][][.8]{$m$}
		\psfrag{3}[][][.4]{$m{+}1$}
		\psfrag{-}[][][.4]{$m{-}1$}
	\centering\includegraphics[width=9cm]{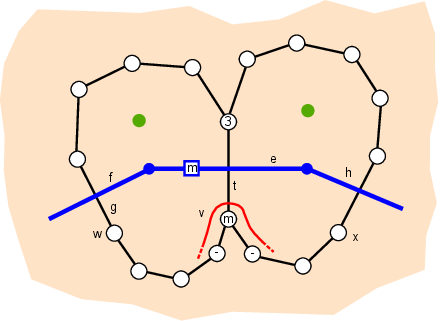}
	\caption{Proof that~$\fu$ is a unicellular map. The part of level~$m$-loop is represented in red.}
	\label{unicel}
\end{figure}

Without loss of generality, let us assume that $v'=v$. Let us furthermore suppose for now that the other extremity of~$\tilde e'$ has label $m-1$. In this case, the level loop is still present and we must have $v''=v$. If the other extremity of~$\tilde e''$ also has label $m-1$, the level loop is still present and we obtain a contradiction. We must thus have that the other extremity of~$\tilde e''$ also has label $m+1$ and as a consequence, $e''$ has label~$m$. We iterate the argument with the blue edges succeeding~$e''$ in the cycle and obtain that every edge of the cycle has label~$m$, a contradiction with our hypothesis that an extremity of~$\tilde e'$ has label $m-1$. In the end, we obtain that~$e'$ has label~$m$ and, iterating the argument, we finally obtain that the blue cycle turns around a single vertex that is a local minimum for the labels. This vertex is thus necessarily~$v^\ooo$.

We can then retract to~$0$ each connected component of the complement of~$\fu$ along the corresponding connected component of the blue graph. As a result, we obtain that~$\fu$ is a map and we conclude that it is unicellular by Euler's characteristic formula, noting that it has $v(\m)-1+f(\m)$ vertices and $e(\m)$ edges by Proposition~\ref{faces}, where we denoted by $v(\m)$, $f(\m)$ and~$e(\m)$ the number of vertices, faces and edges of~$\m$.
\end{pre}

We are now ready to prove that our orientation process, and thus~$\Phi$, are well defined.

\begin{pre}[Proof of proposition~\ref{oriterm}]
Let us consider the set~$\cC$ of oriented corners~$\vc$ of~$\m$ such that $\lab(\varphi(\vc))<\lab(\vc)$. The proposition follows from the claim that $\cC=\{\upd^i(\vrho),\opp(\si(\upd^i(\vrho)))\,:\,i\ge 0\}$. Indeed, after updating the oriented corner~$\vc$, the level loop issued from~$\si(\vc)$ is oriented, as well as that issued from~$\upd(\vc)$, as it is the same as that issued from~$\varphi^{-k}\circ\si^j(\vc)$, with the notation of Figure~\ref{update}. As every level loop except the one at level~$0$ is issued from at least one corner of~$\cC$, the claim entails that every loop is oriented by the process, which thus terminates. (Note that $\pr^{\lab(\vrho)}(\vrho)$ is incident to~$v^\ooo$ so that the level~$0$-loop is always oriented at the beginning of the process.) 

First, observe that $\vrho\in\cC$ (recall that we have chosen between~$\m$ and~$\bar\m$ the map with this property) and that~$\cC$ is stable by~$\upd$ and by~$\opp\circ\si$. Now our orientation process orients several loops, possibly not all the loops. Consider the labeled unicellular map~$(\fu,\lab)$ obtained from the construction of Section~\ref{secmu}, using our orientation process for the loops that our orientation process orients and an arbitrary orientation for the loops that are not oriented by our process. 

From Proposition~\ref{faces}.\ref{faces_iv}, we see that there is a one-to-one correspondence between the oriented corners of~$\cC$ and the white oriented corners of~$\fu$ and, if $\vc\in\cC$ corresponds to the oriented corner~$\va$ of~$\fu$, then $\opp(\si(\vc))$ corresponds to~$\opp(\va)$. Next, with the notation of Figure~\ref{update}, note that all the~$j$ black edges between~$\vc$ and~$\si^j(\vc)$ belong to the same white corner of~$\fu$ and the same goes for the~$l$ black edges between $\varphi^{-k}\circ\si^j(\vc)$ and $\upd(\vc)$. Furthermore, as there are no stops between~$\si^j(\vc)$ and $\varphi^{-k}\circ\si^j(\vc)$ along the contour of the face of~$\m$, by construction, the stops at~$\si^j(\vc)$ and $\varphi^{-k}\circ\si^j(\vc)$ will be linked to the corresponding green vertex during Step~\ref{step3} without any green edges in between. As a result, the oriented corner of~$\fu$ corresponding to~$\upd(\vc)$ is~$\varphi^2(\va)$. Consequently, the sequence of oriented corners of~$\fu$ corresponding to $\upd^i(\vrho)$, $i\ge 0$ are the white corners of the unique face (Lemma~\ref{lemul}) of~$\fu$, arranged in the root-induced contour order. Consequently, every oriented corner of~$\fu$ is either of this form or its opposite is. The claim follows.
\end{pre}

We may now state our main theorem.

\begin{thm}\label{bijthm}
The mappings $\Phi:\cB^\ooo \to \cU\times\{+,-\}$ and $\Psi:\cU\times\{+,-\} \to \cB^\ooo$ are bijections, which are inverse one from another.
\end{thm}

Using Proposition~\ref{faces}, we obtain the following specialization of Theorem~\ref{bijthm}. For an integer finite sequence $\alpha=(\alpha_1, \ldots, \alpha_n)$, we denote by $\cB^\ooo_{\alpha}$ the set of pointed bipartite maps of~$\cS$ with~$n$ faces marked~$1$, $2$, \ldots, $n$ such that, for $1\le i\le n$, the face marked~$i$ has degree~$2\alpha_i$. We also denote by~$\cU_\alpha$ the set of well-labeled unicellular mobiles with~$n$ green vertices marked~$1$, $2$, \ldots, $n$ and such that, for $1\le i\le n$, the green vertex marked~$i$ has degree~$\alpha_i$.
\nomenclature[9ab]{$\cB^\ooo_{\alpha}$}{set of pointed bipartite maps of~$\cS$ with face degrees $2\alpha_1$, \dots, $2\alpha_n$}%
\nomenclature[9bb]{$\cU_\alpha$}{set of well-labeled unicellular mobiles of~$\cS$ with green vertex degrees~$\alpha_1$, \dots, $	\alpha_n$}%

\begin{corol}\label{bijcor}
The restrictions of~$\Phi$ to $\cB^\ooo_{\alpha}$ and of~$\Psi$ to $\cU_\alpha\times\{+,-\}$ are bijections, which are inverse one from another.
\end{corol}

\begin{pre}[Proof of Theorem~\ref{bijthm}]
It is sufficient to show that
\begin{enumerate}[label=(\textit{\alph*})]
	\item If $(\fu,\lab)\in\cU$, $\eps\in\{+,-\}$ and $(\m,v^\ooo)=\Psi((\fu,\lab),\eps)$, then $\Phi(\m,v^\ooo)=((\fu,\lab),\eps)$.\label{bij_a}
	\item If $(\m,v^\ooo)\in \cB^\ooo$ and $((\fu,\lab),\eps)=\Phi(\m,v^\ooo)$, then $(\fu,\lab)\in\cU$ and $\Psi((\fu,\lab),\eps)=(\m,v^\ooo)$.\label{bij_b}
\end{enumerate}

\noindent\textbf{Let us show }\ref{bij_a}\textbf{.}\quad
Let $(\fu,\lab)\in\cU$ and let us denote by $(\m,v^\ooo)\de\Psi((\fu,\lab),-)$. It is enough to show that $\Phi(\m,v^\ooo)=((\fu,\lab),-)$. Indeed, $\Psi((\fu,\lab),+)=(\bar\m,v^\ooo)$ and the previous fact immediately implies that $\Phi(\bar\m,v^\ooo)=((\fu,\lab),+)$. 

As $(\m,v^\ooo)\in \cB^\ooo$ (Proposition~\ref{propmu}), we can apply to it the construction of Section~\ref{secmu}. First, notice that the labels of~$V(\m)$ are given by~$\lab$. Indeed, $\lab(v^\ooo)=0$ and, if $v\neq v^\ooo$, then, as the variation of~$\lab$ along any black edge is~$1$, every black path from~$v^\ooo$ to~$v$ has length at least $\lab(v)$. Moreover, starting from an arbitrary corner of~$\fu$ incident to~$v$ and following the black edges linking it to its iterate successors provides a path from~$v$ to $v^\ooo$ of length~$\lab(v)$. In particular, we furthermore immediately obtain that the root vertex of~$\m$ has a larger label than the other extremity of the root edge of~$\m$ (as on Figure~\ref{rootm}), which entails that the second coordinate of $\Phi(\m,v^\ooo)$ will be set to~$-$.

From the proof of Proposition~\ref{propmu}, we see that each face of~$\m$ corresponds to a unique green vertex of~$\fu$. As the rooting conventions (Figures~\ref{rootm} and~\ref{rootu}) clearly correspond, it only remains to see that the stops of the loops of~$\m$ are exactly at the locations of the green edges of~$\fu$. 

Let us consider an integer $i\ge 1$ and the corner cycle of a label $i+1$-corner of~$\fu$. We claim that this corner cycle corresponds to a level~$i$-loop of~$\m$. Let us first consider a corner run of this corner cycle. By definition, the labels of all the corners in this run are larger than $i+1$ and the first corners from both sides outside the run are labeled~$i$ or less. Let us denote by~$c_1$ and~$c_2$ the latter corners. We also denote by~$e_1$ (resp.~$e_2$) the green edge linking the vertex incident to~$c_1$ (resp.\ to~$c_2$) to the green vertex at the corresponding extremity of the run; see Figure~\ref{ccll} for visual aid.

\begin{enumerate}[label=(\textit{\roman*})]
\item Let us first assume that the run is trivial, that is, consists of a single green vertex. As in the proof of Proposition~\ref{propmu}, we consider the chain of black edges linking~$c_1$ to~$c_2$, which consists in taking edges linking the iterate successors of~$c_1$ and~$c_2$ until we reach a common corner. This chain of black edges is part of the contour of the face of~$\m$ that corresponds to the unique vertex of the run and the labels along it are all smaller than or equal to~$i$. As a result, if there is a loop at level~$i$ in~$\m$ that crosses~$e_1$, then it follows the chain of black edges and crosses~$e_2$.

\item\label{crnt} Let us now assume that the run is not trivial. As the successor of any corner of the run with label strictly larger than~$i$ belongs to the run or to $\{c_1,c_2\}$, the fact that $(\fu,\lab)$ is well labeled ensures that the range of labels in the run is of the form $\{i+1,i+2,\ldots,k\}$ for some $k\ge i+1$. As a result, there is at least one corner labeled $i+1$ in the run. Moreover, all the corners labeled $i+1$ in the run have the same successor, which is either~$c_1$ or~$c_2$ and must have label~$i$. Up to exchanging roles, let us assume that~$c_1$ is this corner and let $j\de\lab(c_2)\le i$. The same reasoning as in the proof of Proposition~\ref{propmu} may be conducted with~$c_1$ and~$c_2$, and we obtain a chain of either $i-j$ or $i-j+2$ black edges (depending on whether~$c_2$ is an iterate successor of~$c_1$ or not) that links~$c_1$ to~$c_2$. Now there is level~$i$-loop of~$\m$ that crosses~$e_1$; it then crosses all the black edges linking the label $i+1$-corners of the run to~$c_1$ and finally follows the chain of black edges until it crosses~$e_2$.
\end{enumerate}
Crossing a green edge for a level loop amounts to applying~$\sigma$ in the algorithm that defines corner cycles, so that the considered corner cycle indeed corresponds to a level~$i$-loop as claimed. Conversely, every level loop at positive level that is not at maximal level in a face crosses at least one black edge linking a corner~$c$ to its successor. Conducting the above reasoning with the corner cycle of~$c$, we obtain a one-to-one correspondence between corner cycles and level loops that are not maximal level loops (the level~$0$-loop corresponds to the corner cycle of the label~$1$-corners).

\begin{figure}[ht!]
		\psfrag{2}[][][.8]{$i$}
		\psfrag{j}[][][.8]{$j$}
		\psfrag{c}[][][.8]{$e_1$}
		\psfrag{d}[][][.8]{$e_2$}
		\psfrag{3}[][][.5]{$i{+}1$}
		\psfrag{+}[][][.5]{$j{+}1$}
		\psfrag{-}[][][.5]{$j{-}1$}
	\centering\includegraphics[width=.95\linewidth]{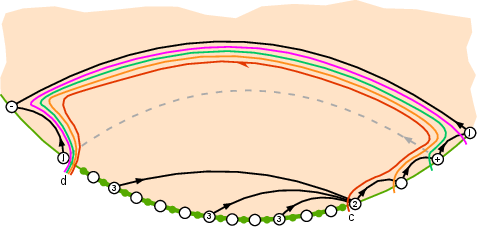}
	\caption{The black edges and the portion of level loop corresponding to a corner run of a label $i+1$-corner, in the case where~$c_2$ is not an iterate successor of~$c_1$. The geodesic orientation of the edges of~$\m$ is represented by full arrowheads. The gray dashed edge represents the case where~$c_2$ is an iterate successor of~$c_1$: in this case, the three bottom-most level loops follow this edge instead of the top one and the pink level loop crosses it. The half arrowhead orients the level~$i$-loop as claimed.}
	\label{ccll}
\end{figure}

As above, we consider a nontrivial corner run and use the notation of~\ref{crnt}. We claim that, once oriented, the corresponding portion of level loop visits the face of~$\fu$ from~$e_1$ to~$e_2$ (see Figure~\ref{ccll}). As the labels along the chain of black edges are all strictly smaller than~$i$ except that of~$c_1$ and possibly that of~$c_2$, the claim entails that the portion of loop might have stops only at the locations of~$e_1$ and~$e_2$. More precisely, it has a stop at~$e_2$ if and only if $j=i$. Moreover, by the same argument, for a trivial run, the corresponding portion of level~$i$-loop may also have stops only at the locations of the green edges and it has one at the exiting green edge if and only if the corresponding label is~$i$. As the entering green edge of the considered portion of loop is the exiting green edge of the previous portion of the same loop, the claim yields that the stops of the loops of~$\m$ are exactly at the locations of the green edges of~$\fu$, as desired.

In order to show the claim, first observe that the coherent orientation process of the corners of~$\fu$ may be restated as follows. Consider all the corner cycles that contain the green corner~$\varphi(\vtau)$ (where~$\vtau$ is the root of~$\fu$) and orient them according to it. Then, iteratively, consider all the corner cycles that have not been oriented and that contain the next green corner in the order $\varphi(\vtau)$, $\varphi^3(\vtau)$, $\varphi^5(\vtau)$, \ldots{} and orient them according to it. Stop when all the green corners have been considered.

The first green corner~$\varphi(\vtau)$ needs a special treatment as it serves to orient all the corner cycles that contain it, in contrast with the other green corners, which may only serve to orient the corner cycles such that the green corner is an extremity of one of its runs. This special treatment corresponds to the initialization step of the orientation process of the level loops. The corner cycles of the iterate successors of the root~$\vtau$, $\suc(\vtau)$, \ldots, $\suc^{\lab(\vtau)-1}(\vtau)$ are exactly the ones that contain the root corner; their coherent orientation is thus the same as their root-induced orientation. For $i\in\{1,\ldots,\lab(\vtau)\}$, the loop issued from $\pr^{i}(\vrho)$ (where~$\vrho$ is the root of~$\m$) corresponds to the corner cycle of $\suc^{i-1}(\vtau)$ and the initialization step of the orientation process of the loops orients them as desired.

As long as the orientation of the loops matches the orientation of the corner cycles, the green edges and the exploration of the green corners of $\Phi(\m, v^\ooo)$ match the green edges and the exploration of the green corners of~$\fu$. Let us now consider a green corner. Then all the level loops corresponding to the corner cycles containing the green corner as a run extremity are considered in the update step of the orientation process of the level loops (as depicted on Figure~\ref{update}). Furthermore, if not oriented, the level loops are oriented during this step in the desired orientation. The claim follows.

\bigskip
\noindent\textbf{We now turn to }\ref{bij_b}\textbf{.}\quad It is sufficient to consider pointed bipartite maps such that the distance from the distinguished point to the root vertex is larger than that to the other extremity of the root edge. Let $(\m,v^\ooo)\in\cB^\ooo$ be such a map and $((\fu,\lab),-)=\Phi(\m,v^\ooo)$.

We have already established that $(\fu,\lab)$ is a labeled unicellular mobile (Lemma~\ref{lemul}) and that there is a unique outgoing edge of~$\m$ for each corner of~$\fu$ (Proposition~\ref{faces}.\ref{faces_iv}). We claim that, for each corner of~$\fu$, the corresponding outgoing edge of~$\m$ links it to its successor. This claim entails that $(\fu,\lab)$ is well labeled as, by definition of~$\lab$, the head label of an edge is one less than its tail label. Furthermore, the construction of Section~\ref{secum} precisely consists in linking with black edges the corners of~$\fu$ to their successors, so that $\Psi((\fu,\lab),-)=(\m,v^\ooo)$.

The claim is obvious for corners with label~$1$. Let us thus fix $i\ge 1$ and consider a level~$i$-loop of~$\m$; see Figure~\ref{mvulmv}. As it has at least one stop, it has to cross a green edge. We consider a portion of this level loop completely included in the face of~$\fu$: it starts from a green edge~$e_1$, follows black edges, possibly crossing \edg{$i$}{$i$+$1$} black edges along the way and finish on a green edge~$e_2$. The chain of black edges along which the portion of loop evolves splits the face of~$\fu$ into two connected components; let us focus on the one that contains the portion of loop. The boundary of this component consists of the chain of black edges, $e_1$, $e_2$, and a chain of green edges whose extremities are the green vertices incident to~$e_1$ and~$e_2$. Let us denote by~$c_1$ and~$c_2$ the white corners of~$\fu$ through which the portion of loop enters and exits the face of~$\fu$. By definition of a level loop, the labels along the black chain are all smaller than or equal to~$i$. More precisely, because of the definition of stops, we have $\lab(c_1)\le i$, $\lab(c_2)\le i$, and the other labels along the black chain are all strictly smaller than~$i$.

If the green chain contains at least a white corner, then the portion of loop has to cross an \edg{$i$}{$i$+$1$} black edge, as otherwise the map would not be connected or a face would contain two green vertices. By definition of stops, all the \edg{$i$}{$i$+$1$} black edges crossed by the portion of loop must end in~$c_1$. We know that there is an outgoing edge of~$\m$ in every corner of~$\fu$; subsequently following these outgoing edges from any corner of~$\fu$ yields a path with decreasing labels linking the corner to~$v^\ooo$. Starting from a corner on the green chain, this path must intersect the black chain and this may only happen at~$c_1$. As a result, the range of labels on the green chain is of the form $\{i+1,i+2,\ldots,k\}$ for some $k\ge i+1$. The green chain is thus a corner run of all the label~$i+1$-corners of the chain. Moreover, provided that the coherent orientation of the run is from~$e_2$ to~$e_1$, then~$c_1$ is the successor of all these corners.

\begin{figure}[ht!]
		\psfrag{i}[][][.8]{$i$}
		\psfrag{j}[][][.45]{$\le i$}
		\psfrag{z}[][][.8]{$c$}
		\psfrag{c}[][][.8]{$e_1$}
		\psfrag{d}[][][.8]{$e_2$}
		\psfrag{-}[][][.45]{$<i$}
	\centering\includegraphics[width=13cm]{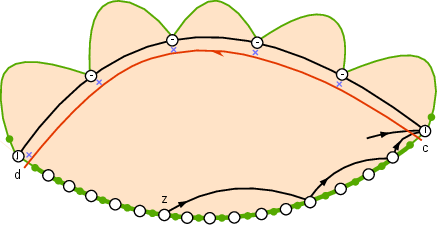}
	\caption{A portion of a level~$i$-loop completely included in the face of~$\fu$ delimits a corner run of label~$i+1$-corners.
	The purple crosses mean that there cannot be black edges. From a corner~$c$ on the green chain, following the successive outgoing edges yields a geodesic to~$v^\ooo$: it has to cross the black chain and this may only happen at~$c_1$.}
	\label{mvulmv}
\end{figure}

We obtain that the level loops of~$\m$ correspond to the corner cycles of~$(\fu,\lab)$. It only remains to see that the level loops enter the face of~$\fu$ from the end of the corresponding coherently oriented corner runs. Let us consider a corner cycle and the first green corner~$c$ in the root-induced contour order it contains. It suffices to see that the corner run of~$c$ is coherently oriented in accordance with the corresponding portion of level loop, that is, from the exiting green edge to the entering green edge. On the one hand, the corner run is oriented by the root-induced orientation of~$c$. On the other hand, observe that our orientation process for the level loops explores the green corners in the root-induced contour order and orients the newly explored loops in the orientation opposite to that of the root~$\vtau$ of~$\fu$. If the green corner~$c$ corresponds to~$\varphi(\vtau)$, then the corresponding level loop is one issued from an iterate parent of the root of~$\m$ and its orientation is opposite to that of~$\vtau$. Otherwise, the loop corresponding to the considered corner cycle will be oriented when we explore~$c$, in the orientation opposite to the root-induced orientation of~$c$, as desired. The claim follows.
\end{pre}

\section{General maps}\label{secgen}

\subsection{The mappings}

We now relax the hypothesis that the map is bipartite. We will slightly modify it in order to be able to apply our bijection from Section~\ref{secmu}. We denote by~$\M^\ooo$ the set of pointed maps of~$\cS$ and by~$\M^\ooo_{\textrm{eq}}$ its subset consisting of pointed maps such that both extremities of the root edge are at the same distance from the distinguished vertex.
\nomenclature[Aa]{$\M^\ooo$}{set of pointed maps of~$\cS$}%
\nomenclature[Ab]{$\M^\ooo_{\textrm{eq}}$}{set of pointed maps of~$\cS$ whose root edge is equilabeled}%

\begin{figure}[ht!]
		\psfrag{v}[][]{$v^\ooo$}
		\psfrag{a}[][][.8]{$a$}
		\psfrag{b}[][][.8]{$b$}
		\psfrag{c}[][][.8]{$c$}
		\psfrag{d}[][][.8]{$d$}
		\psfrag{e}[][][.8]{$e$}
		\psfrag{0}[][][.8]{$0$}
		\psfrag{1}[][][.8]{$1$}
		\psfrag{2}[][][.8]{$2$}
		\psfrag{3}[][][.8]{$3$}
		\psfrag{4}[][][.8]{$4$}
	\centering\includegraphics[width=.95\linewidth]{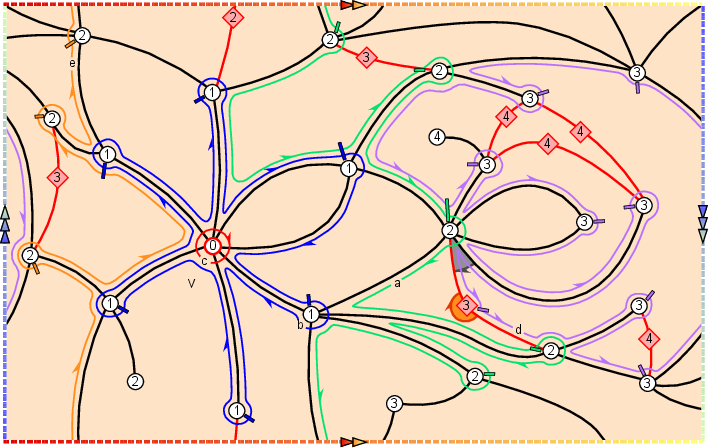}
	\caption{The level loops for a general map. The equilabeled edges are represented in red and the added vertices are represented by red squares. As above, the maximal level loops are not represented. Beware that, in this example, the root of~$\tilde\m$ (in gray as usual) is such that the root vertex is closer to~$v^\ooo$ than the other extremity of the root edge; one thus has to consider the root-flipped map whose root is represented in orange.}
	\label{bij_gen1}
\end{figure}

We consider a general pointed map $(\m,v^\ooo)\in\M^\ooo$ and define the function $\lab:V(\m)\to\Z_+$ by~\eqref{defl} as before. There are now two kinds of edges: an edge will be called \emph{equilabeled} if its extremities both have the same label. We then enlarge the map~$\m$ by adding in the middle of each equilabeled edge an extra vertex splitting the edge into two new edges. We denote by~$\tilde\m$ this enlarged map and we assign to each added vertex the common label of its two neighbors plus~$1$. This extends the definition of~$\lab$ to $V(\tilde\m)$ and, clearly, for $v\in V(\tilde\m)$, one has $\lab(v)= d_{\tilde\m}(v^\ooo,v)$; see Figure~\ref{bij_gen1}.

\begin{figure}[ht!]
		\psfrag{v}[][]{$v^\ooo$}
		\psfrag{0}[][][.8]{$0$}
		\psfrag{1}[][][.8]{$1$}
		\psfrag{2}[][][.8]{$2$}
		\psfrag{3}[][][.8]{$3$}
		\psfrag{4}[][][.8]{$4$}
	\centering\includegraphics[width=.95\linewidth]{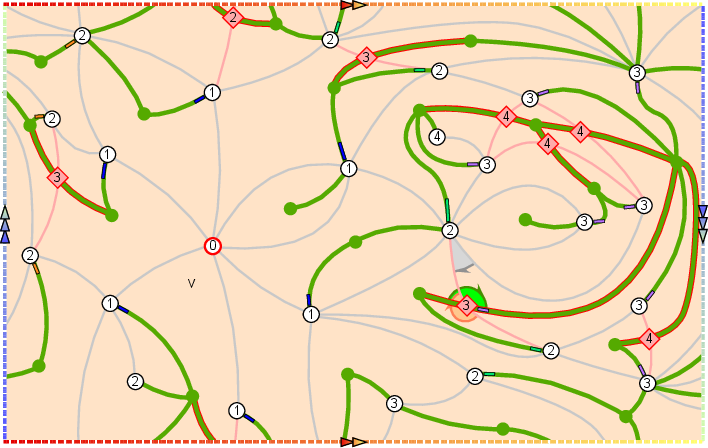}
	\caption{The bijection for a general map. The flagged edges of the mobile have been highlighted and their labels are represented by the red squares. In this example, the root edge is equilabeled; the root of the mobile is thus a flagged edge given with a side and local orientation.}
	\label{bij_gen2}
\end{figure}

The map~$\tilde\m$ is bipartite: we apply to it the construction of Section~\ref{secmu} and set $((\tilde\fu,\lab),\eps)\de \Phi(\tilde\m,v^\ooo)$. Note that we necessarily have $\eps=+$ in the case $(\m,v^\ooo)\in\M^\ooo_{\textrm{eq}}$. We slightly modify the encoding map as follows. Every vertex of $V(\tilde\m)\setminus V(\m)$ is by design of degree~$2$ and there are stops at both corners incident to it, as their labels are local maximums along the boundaries of the incident faces. As a result, it also has degree~$2$ in~$\tilde\fu$; we suppress it and merge the two incident edges into a single edge. We call such a resulting edge a \emph{flagged edge} and we assign to it the label of the suppressed vertex. We denote by $(\fu,\lab)$ the resulting map, the function~$\lab$ being defined on a subset of the vertices and edges of~$\fu$. Finally, if the root edge of~$\m$ is equilabeled, then the root vertex of~$\tilde\fu$ is one of the added vertices of~$\tilde\m$. In this case, we transgress our usual definition of root and declare the root of~$\fu$ to be the edge resulting from the suppression of the root vertex of~$\tilde\fu$, together with the side and local orientation induced by the root of~$\tilde\fu$. Such a map will be called \emph{edge-rooted} in what follows; see Figure~\ref{bij_gen2}.

We extend the definition of~$\Phi$ by setting $\Phi(\m,v^\ooo)\de((\fu,\lab),\eps)$. Plainly, if~$\m$ is bipartite, then $\tilde\m=\m$ and $((\tilde\fu,\lab),\eps)=((\fu,\lab),\eps)$ so that the definition of~$\Phi$ is consistent with the previous definition from Section~\ref{secmu}.

\begin{defi}[Labeled generalized unicellular mobile]
A \emph{labeled generalized unicellular mobile} is a pair~$(\fu,\lab)$ such that
\begin{itemize}
	\item $\fu$ is a rooted or edge-rooted one-face map of~$\cS$ whose vertex set is partitioned into $V_\oov(\fu)\sqcup V_\circ(\fu)$ in such a way that every edge has at least one extremity in $V_\oov(\fu)\,$;
	\item if~$\fu$ is rooted, its root vertex lies in~$V_\circ(\fu)$; if~$\fu$ is edge-rooted, its root edge belongs to the set $E_\sqr(\fu)$ of edges linking two vertices of $V_\oov(\fu)\,$;
	\item $\lab:V_\circ(\fu)\sqcup E_{\sqr}(\fu)\to\N$ is a function with minimum~$1$.
\end{itemize}
\end{defi}
\nomenclature[Ag]{$E_\sqr(\fu)$}{set of flagged edges of~$\fu$}%

The edges of $E_\sqr(\fu)$ are called \emph{flagged edges}. We define the following \emph{flag splitting operation} on the set of labeled generalized unicellular mobiles. Given a labeled generalized unicellular mobile, we split each flagged edge into two edges by adding a vertex in the middle and we assign to the added vertex the label of the flagged edge. If the original map is edge-rooted, we root the resulting map at the corner incident to the vertex added on the root edge that corresponds to the distinguished side of the root edge, oriented in accordance with the root orientation of the original map. We obtain by this operation a labeled unicellular mobile.

\begin{defi}[Well-labeled generalized unicellular mobile]
A labeled generalized unicellular mobile is called a \emph{well-labeled generalized unicellular mobile} if the flag splitting operation makes it a well-labeled unicellular mobile and the corners incident to vertices added during the flag splitting operation are not successors of any corners.
\end{defi}

Let $(\fu,\lab)$ be a well-labeled generalized unicellular mobile, let $(\tilde\fu,\lab)$ denote the well-labeled unicellular mobile obtained by the flag splitting operation and let $(\tilde\m,v^\ooo)\de \Psi((\tilde\fu,\lab),+)$. Let $v\in V(\tilde\fu)\setminus V(\fu)$. It has degree~$2$ in~$\tilde\fu$ and, by definition, the two incident corners are not successors of any corners, so that it also has degree~$2$ in~$\tilde\m$. As a result, we may define the map~$\m$ by suppressing from~$\tilde\m$ each vertex of $V(\tilde\fu)\setminus V(\fu)$ and by merging the two incident edges into a single edge. Note that neither~$v^\ooo$ nor the root vertex of~$\tilde\m$ are suppressed, so that the pointed map $(\m,v^\ooo)$ is well defined.

We extend the definition of~$\Psi$ by setting $\Psi((\fu,\lab),+)\de(\m,v^\ooo)$ and, whenever $(\fu,\lab)$ is not edge-rooted, $\Psi((\fu,\lab),-)\de(\bar\m,v^\ooo)$. We denote by~$\calG$ the set of well-labeled generalized unicellular mobiles and by $\calG_{\textrm{er}}\subset\calG$ the subset of edge-rooted ones. If $(\fu,\lab)\in\calG$, the \emph{generalized degree} of a vertex in $V_\oov(\fu)$ is defined as the number of incident nonflagged edges plus half the number of incident flagged edges (counted with multiplicity). We readily obtain the following theorem.
\nomenclature[Ad]{$\calG$}{set of well-labeled generalized unicellular mobiles}%
\nomenclature[Ae]{$\calG_{\textrm{er}}$}{set of edge-rooted well-labeled generalized unicellular mobiles}%

\begin{thm}\label{bijthmgen}
The extended mappings~$\Phi$ and~$\Psi$ induce bijections between~$\M^\ooo_{\textrm{eq}}$ and $\calG_{\textrm{er}}\times\{+\}$ on the one hand, and between~$\M^\ooo\setminus\M^\ooo_{\textrm{eq}}$ and $\calG\setminus\calG_{\textrm{er}}\times\{+,-\}$ on the other hand, and these bijections are inverse one from another.

Moreover, if $(\m,v^\ooo)\in\M^\ooo$ and $((\fu,\lab),\eps)=\Phi(\m,v^\ooo)$, then
\begin{enumerate}[label=(\textit{\roman*})]
	\item $V(\m)=V_\circ(\fu)\sqcup \{v^\ooo\}$ and, for $v\in V_\circ(\fu)$, $\lab(v)=d_\m(v,v^\ooo)\,$;
	\item the faces of~$\m$ correspond to $V_\oov(\fu)\,$; moreover, the degree of a face of~$\m$ is twice the generalized degree of the corresponding vertex in $V_\oov(\fu)$;
	\item the maps~$\m$ and~$\fu$ have the same number of edges.
\end{enumerate}
\end{thm}

We may specialize this theorem by prescribing the face degrees as we did in Corollary~\ref{bijcor}. For a finite integer sequence $\alpha=(\alpha_1, \ldots, \alpha_n)$, we denote by $\M^\ooo_{\alpha}$ the set of pointed maps of~$\cS$ with~$n$ faces marked~$1$, $2$, \ldots, $n$ such that, for $1\le i\le n$, the face marked~$i$ has degree~$\alpha_i$. We also denote by~$\calG_\alpha$ the set of well-labeled generalized unicellular mobiles with~$n$ green vertices marked~$1$, $2$, \ldots, $n$ and such that, for $1\le i\le n$, the green vertex marked~$i$ has generalized degree~$\alpha_i/2$.
\nomenclature[Ac]{$\M^\ooo_{\alpha}$}{set of pointed maps of~$\cS$ with face degrees~$\alpha_1$, $\alpha_2$, \ldots, $\alpha_n$}%
\nomenclature[Af]{$\calG_\alpha$}{set of elements of $\calG$ with green vertex generalized degrees~$\alpha_1/2$, \ldots, $\alpha_n/2$}%

\begin{corol}
The restriction of $\Phi$ to $\M^\ooo_{\alpha}$ realizes a bijection between~$\M^\ooo_{\alpha}$ and $(\calG_\alpha\times\{+,-\})\setminus(\calG_{\textrm{er}}\times\{-\})$.
\end{corol}

\subsection{Application to triangulations}\label{sectri}

\subsubsection{The bijection}

For quadrangulations, the bijections are more convenient: in particular, the labels on the encoding objects only satisfy \emph{local} constraints instead of global ones as it is the case in general (see Appendix~\ref{secquad}). As explained at the end of Section~\ref{seccoh}, the orientation of the corner cycles does not intervene in the definition of the encoding objects. In the case of triangulations, the same argument holds and we can thus give a simpler characterization of the encoding objects and derive enumeration results.

\begin{figure}[ht!]
		\psfrag{-}[][][.5]{$i{-}1$}
		\psfrag{i}[][][.8]{$i$}
		\psfrag{+}[][][.5]{$i{+}1$}
	\centering\includegraphics[width=.95\linewidth]{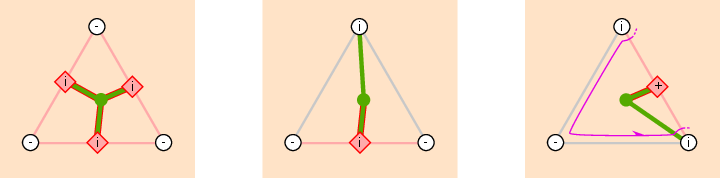}
	\caption{The three possible types of face and the corresponding green vertices.}
	\label{labelstri}
\end{figure}

Let us look at the restriction of~$\Phi$ to the set~$\T^\ooo$ of pointed triangulations. By Theorem~\ref{bijthmgen}, every green vertex of a corresponding element of~$\calG$ is either exactly incident to~$3$ flagged edges or exactly incident to one nonflagged edge and one flagged edge. Moreover, the labels around a green vertex can only be of the three types depicted on Figure~\ref{labelstri}.
\nomenclature[Ba]{$\T^\ooo$}{set of pointed triangulations of~$\cS$}%

Let us denote by~$\EE$ the set of labeled generalized unicellular mobiles with green vertices of generalized degree~$3/2$ and such that the labels around any green vertex are of a type shown on Figure~\ref{labelstri}. We consider an element of~$\EE$ and apply to it the flag splitting operation. Plainly, the label variations between adjacent corners (corners separated by exactly one green corner) belong to $\{-1,0,1\}$. As in the case of quadrangulations, this entails that the requirements of Definition~\ref{wlum} are fulfilled, no matter how the corners are oriented. Moreover, the conditions on the labels entail that the first vertices outside any corner run are original white vertices (in the sense that they were not added during the flag splitting operation) so that the successors will always be incident to original vertices. Denoting by~$\EE_{\textrm{er}}$ the set of edge-rooted mobiles from~$\EE$, we obtain the following corollary.
\nomenclature[Bb]{$\EE$}{set of labeled generalized unicellular mobiles with green vertices as on Figure~\ref{labelstri}}%
\nomenclature[Bc]{$\EE_{\textrm{er}}$}{set of edge-rooted mobiles from~$\EE$}%

\begin{corol}\label{cortri}
The restriction of~$\Phi$ to $\T^\ooo$ realizes a bijection between~$\T^\ooo$ and $(\EE\times\{+,-\})\setminus(\EE_{\textrm{er}}\times\{-\})$. Moreover, the corresponding objects have the same number of edges.
\end{corol}

\subsubsection{Generating functions}

From now on, it will be convenient to slightly modify the bijection of Corollary~\ref{cortri} by subtracting $1/2$ from all the labels of the flagged edges. This will bring some symmetry and make the computation easier. Note that this also makes sense from a metric point of view, as the labels now represent the distances to the distinguished vertex in the metric space obtained from the triangulation by replacing each edge by a unit length segment. We furthermore consider from now on labeled generalized unicellular mobiles up to addition of an integer constant to all the labels. We will decompose these mobiles into elementary pieces in order to count them. We use the framework of~\cite{chapuy07brm,ChDo15bij}.

\begin{figure}[ht!]
		\psfrag{=}[][][.7]{$=$}
		\psfrag{+}[][][.7]{$+$}
		\psfrag{i}[][][.7]{$i$}
		\psfrag{j}[l][l][.7]{$i\pm\frac12$}
	\centering\includegraphics[width=.95\linewidth]{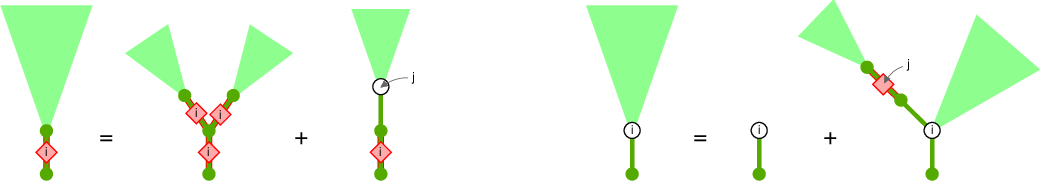}
	\caption{Elementary decomposition of $F$-mobiles and $N$-mobiles.}
	\label{FNdecomp}
\end{figure}

We first consider nonrooted plane labeled generalized unicellular mobiles with one green vertex incident to exactly one edge and such that every other green vertex is incident to exactly
\begin{itemize}
	\item either~$3$ flagged edges with the same label
	\item or one flagged edge with some label $i\in\Z+\frac12$ and one nonflagged edge incident to a white vertex with label $i\pm \frac12$.
\end{itemize}
We introduce the generating function~$F$ (resp.~$N$) counting with weight~$t$ per edge (flagged or nonflagged) such objects with the extra condition that the isolated green vertex is incident to a flagged (resp.\ nonflagged) edge. The decomposition depicted on Figure~\ref{FNdecomp} yields
\begin{equation}\label{FNeq}
F=tF^2+2tN\qquad\text{ and }\qquad N=t+2tNF.
\end{equation}
\nomenclature[Ca]{$F$}{generating function of $F$-mobiles}%
\nomenclature[Cb]{$N$}{generating function of $N$-mobiles}%

Let us call \emph{$F$-mobile} and \emph{$N$-mobile} the objects counted by~$F$ and $N$. When we fix a representative of the labels of an $F$-mobile (resp.\ an $N$-mobile), we call label of the mobile the label of the flagged edge (resp.\ of the white vertex incident to the nonflagged edge) that is incident to the isolated green vertex.

For $i$, $j\in\Z$, we introduce the generating function $C_{i,j}$ of nonrooted plane generalized unicellular mobiles whose green vertices satisfy the previous itemized conditions, with two \textbf{distinct} distinguished white vertices labeled~$i$ and~$j$. We extend the definition to~$i\in\Z+\frac12$ by replacing the distinguished white vertex labeled~$i$ by a green vertex incident to exactly one flagged edge labeled~$i$. Similarly, we extend the definition to every $i$, $j\in\Z/2$.

In order to compute $C_{i,j}$, we start by considering the following objects. We link the isolated green vertices of two $N$-mobiles by a chain of flagged edges and we graft an $F$-mobile onto every green vertex of the chain that is incident to two flagged edges. Every $F$-mobile can be grafted either on one side or on the other side of the chain. We furthermore impose that the labels of the resulting plane mobile satisfy the previous conditions. If we let~$i$ and~$j$ be the labels of the $N$-mobiles, the only possibilities are $j\in\{i-1,i,i+1\}$. Moreover, all the flagged edges of the chain and all the grafted $F$-mobiles must have the same label, which is $i+1/2$ if $j=i+1$, $i-1/2$ if $j=i-1$, or $i\pm 1/2$ if $j=i$. As a result, the generating function of these objects is $tN^2/(1-2tF)=N^3$ if $j=i\pm 1$ and $2N^3$ if $j=i$.

We then consider Motzkin words, that is, finite sequences of $-1$, $0$ or $+1$. We count them with a weight~$Z\de N^3$ per $\pm 1$ and $2Z$ per~$0$. Let us call \emph{increment} of a Motzkin word $w_1 w_2\ldots w_n$ the integer $\sum_{k=1}^n w_k$. We denote by~$U$ the generating function of words $w_1 w_2\ldots w_n$ with increment~$-1$ and such that $\sum_{k=1}^l w_k\ge 0$ for $1\le l<n$, and by~$B$ the generating function of words with increment~$0$. Considering whether the first letter is $-1$, $0$ or $+1$, we obtain
$$U=Z+2Z\,U+Z\,U^2\qquad\text{ and }\qquad B=1+Z\,UB+2Z\,B+Z\,UB$$
so that
\begin{equation}\label{ZBU}
Z=\frac{U}{(1+U)^2}\,,\qquad B=\frac{1+U}{1-U}\,,\qquad U=\frac{1-\sqrt{1-4Z}}{2Z}-1.
\end{equation}
The generating function of Motzkin words with increment~$\ell\in\Z$ is then equal to $B U^{|\ell|}$.

We may now explicit $C_{i,j}$. Let us first suppose that $i$, $j\in\Z$. The differences between the labels of subsequent white vertices on the path linking the two distinguished vertices form a nonempty Motzkin word. Decomposition at these vertices yields $C_{i,j}=BU^{|j-i|}-\un{\{i=j\}}$.
Now, if $i\in\Z+\frac12$, $j\in\Z$, we have $C_{i,j}=N^2(BU^{|j-i-\frac12|}+BU^{|j-i+\frac12|})=B\sqrt{N}U^{|j-i|}$. Similar computations show that, for~$i$, $j\in\Z/2$, we have
\begin{equation}\label{Cij}
C_{i,j}=BN^{\frac12|\{i,j\}\cap(\Z+\frac12)|}\,U^{|j-i|}-\un{\{i=j\in\Z\}}.
\end{equation}

\subsubsection{The sphere}

We may now prove Propositions~\ref{enum} and~\ref{gfun}. We start with the particular case of the sphere, which is somehow degenerate and has to be treated separately. It is relatively easy and has probably already been done by this method but, as we did not find it in the literature, we do it here for future reference.

The series we consider are really series in $x\de t^3$. Setting $\si\de tF$ and using~\eqref{FNeq}, we obtain that
\begin{equation}\label{sigmax}
tF=\si\,,\qquad\frac{N}{t}=\frac{1}{1-2\si}\,,\qquad x=\frac12 \si(1-\si)(1-2\si).
\end{equation}
%
%

\begin{pre}[Proof of Propositions~\ref{enum} and~\ref{gfun} when~$\cS$ is the sphere]
The generating function of rooted encoding mobiles is given by $N/t-1$ and the generating function of edge-rooted encoding mobiles is given by $F^2/t$. By Corollary~\ref{cortri} and~\eqref{FNeq}, the generating function of \textbf{pointed} plane triangulations is thus
$$\frac{F^2}{t}+2\left(\frac{N}{t}-1\right)=\frac{F}{t^2}-2=\frac{\si}{x}-2$$
The function~$\si$ is an algebraic series in~$x$, which has a dominant singularity at $(12\sqrt{3})^{-1}$ and admits the following Puiseux expansion at its singularity:
\begin{equation*}\label{siPuis}
\si=\frac{1}{2}-\frac{\sqrt{3}}{6}-\frac{\sqrt{2}}{6}\big(1-12\sqrt{3}\,x\big)^{1/2} +\OO\big(1-12\sqrt{3}\,x\big).
\end{equation*}
As a result, the generating function of pointed plane triangulations admits the following Puiseux expansion at its singularity:
$$\frac{\si}{x}-2=6\sqrt{3}-8-2\sqrt{6}\big(1-12\sqrt{3}\,x\big)^{1/2} +\OO\big(1-12\sqrt{3}\,x\big).$$
By classical transfer theorems~\cite[Chapter~VI.3]{flajolet09}, this implies that the number of plane triangulations with $3n$ edges (and thus $2n$ faces and $n+2$ vertices) is asymptotically equivalent to
$$-\frac{1}{n}\,2\sqrt{6}\,\frac{n^{-3/2}}{\Gamma(-1/2)}\big(12\sqrt{3}\big)^n=\frac{\sqrt{6}}{\sqrt{\pi}}\,n^{-5/2}\big(12\sqrt{3}\big)^n$$
(the factor $1/n$ comes from depointing). This is the statement of Proposition~\ref{enum}.

Proposition~\ref{gfun} is obtained by noticing that the generating function counting plane triangulations with weight~$x$ per vertex is given by
$$\int x\left(\frac{\si}{x}-2\right) dx=\frac12 \si^3(1-\si)(1-4\si+2\si^2),$$
as claimed.
\end{pre}

\subsubsection{The projective plane}

We now turn to the projective plane, which is also a degenerate case.

\begin{pre}[Proof of Propositions~\ref{enum} and~\ref{gfun} when~$\cS$ is the projective plane]
A unicellular map of the projective plane consists of a single cycle (whose neighborhood forms a M\"obius strip) on the vertices of which trees are grafted. We consider an encoding mobile $(\fu,\lab)$ and we look at the first edge~$e$ belonging to the cycle that we meet when traveling along the face of~$\fu$ starting from the root. We denote by~$e_-$ the extremity of~$e$ that we met before visiting~$e$ and by~$e_+$ the other extremity.

If~$e$ is nonflagged, then~$e_-$ is necessarily white and~$e_+$ is green. Moreover, $e_-$ is incident to~$e$, another edge of the cycle and two trees, one on each side of the cycle. We cut the cycle at~$e_-$, leaving one of the two trees grafted on each extremity of the chain we obtain. The resulting object is similar to the objects counted by $C_{0,0}$ with the difference that the first $N$-mobile of the chain is rooted or edge-rooted. Noticing that there are as many nonflagged edges as white corners in any generalized unicellular mobile, we see that counting rooted mobiles twice and edge-rooted mobiles once amounts to counting nonrooted mobiles with a weight of double their number of edges. By Corollary~\ref{cortri}, the generating function of pointed triangulations whose corresponding mobile is as considered in this paragraph is thus
$$2t \frac{dN}{dt}\frac{C_{0,0}}{N}.$$

If~$e$ is flagged, two cases may happen. If~$e$ is not the root edge, then the mobile consists of a chain counted by $C_{\frac12,\frac12}$ together with a rooted or edge-rooted $F$-mobile. If~$e$ is the root edge, then the mobile is obtained by gluing together (in a nonorientable way) the extremities of a nontrivial chain counted by $C_{\frac12,\frac12}$. The generating function of pointed triangulations whose corresponding mobile is as considered in this paragraph is thus
$$C_{\frac12,\frac12}\, 2t \frac{dF}{dt}+\frac{C_{\frac12,\frac12}}{t}-1.$$ 

A small computation using~\eqref{ZBU}, \eqref{Cij} and~\eqref{sigmax} yields that the generating function of pointed triangulations on the projective plane is
\begin{equation}\label{ptripplane}
\frac{3}{1-6\si+6\si^2}\left(\frac{1-2\si}{\sqrt{1-6\si+6\si^2}}-1+2\si-2\si^2\right),
\end{equation}
which has a dominant singularity at $(12\sqrt{3})^{-1}$ and admits the following Puiseux expansion:
$$2^{-3/4}\, 3^{5/4} \big(1-12\sqrt{3}\,x\big)^{-3/4}\left(1+\OO\big(\big(1-12\sqrt{3}\,x\big)^{1/2}\big)\right).$$
By transfer theorems, the number of triangulations of the projective plane with~$3n$ edges (and thus~$2n$ faces and $n+1$ vertices) is asymptotically equivalent to
$$\frac{1}{n}\,2^{-3/4}\, 3^{5/4}\,\frac{n^{-1/4}}{\Gamma(3/4)}\big(12\sqrt{3}\big)^n,$$
as desired.
The generating function of triangulations counted with weight~$x$ per vertex is then given by integrating~\eqref{ptripplane} with respect to~$x$.
\end{pre}

\subsubsection{The general case}

We suppose from now on that $h\ge 1$, that is, we exclude the already treated cases of the sphere and the projective plane. The general case is similar to~\cite{chapuy07brm} (see also~\cite{ChDo15bij}) so we only give the main arguments and refer the reader to these references for more detail.

\begin{defi}[Scheme]
A \emph{scheme} is a unicellular map of~$\cS$ with only vertices of degree~$3$ or more. A \emph{normalized scheme} is a pair $(\s,\lab^\star)$, where~$\s$ is a scheme and $\lab^\star:V(\s)\to\Z_+/2$ is such that, if we denote by $\lab_1^\star< \ldots <\lab_k^\star$ the values of its range, then $\lab_1^\star\in\{0,\frac12\}$ and $\lab_{i+1}^\star-\lab_{i}^\star\in\{\frac12,1\}$ for $1\le i<k$.
We denote by $\Sg$ the (finite) set of normalized schemes of~$\cS$.
\end{defi}
\nomenclature[D]{$\Sg$}{set of normalized schemes of~$\cS$}%

We consider an encoding mobile $(\fu,\lab)$. We denote by~$\tilde\s$ its $3$-core, that is, the nonrooted unicellular map obtained by iteratively removing from~$\fu$ all its vertices of degree~$1$ and then replacing the chains of edges linking vertices of degree~$3$ or more by single edges. Note that the vertices of~$\tilde\s$ may be green or white and that it is possible that several edges have both extremities that are white vertices. Note also that green vertices of~$\tilde\s$ necessarily have degree~$3$. We give labels to the vertices of~$\tilde\s$ as follows: to a white vertex, we give the (integer) label of the corresponding vertex in~$\fu$ and to a green vertex, we give the common (noninteger) label of the three flagged edges incident to the corresponding vertex in~$\fu$. Let us denote by $\lab_1< \ldots <\lab_k$ the different values of these labels. We normalize them by replacing them with labels $\lab_1^\star< \ldots <\lab_k^\star$ uniquely defined by $\lab_1^\star\in\{0,\frac12\}$, $\lab_{i+1}^\star-\lab_{i}^\star\in\{\frac12,1\}$ for $1\le i<k$ and $\lab_{i}^\star-\lab_{i}\in\Z$ for $1\le i\le k$. We let~$\s$ be the same map as~$\tilde\s$ but without distinguishing between green and white vertices and rooted as follows. We start from the root of~$\fu$ and we travel along its unique face until we meet a vertex that corresponds to a vertex of~$\s$. The root of~$\s$ is the corner that is visited at this instant, oriented according to the local orientation we were following. The pair $(\s,\lab^\star)$ is then a normalized scheme, which we call \emph{the normalized scheme} of~$(\fu,\lab)$.

At this point, we fix $(\s,\lab^\star)\in\Sg$ and we denote by $\lab_1^\star< \ldots <\lab_k^\star$ the values of the range of~$\lab^\star$ as above. We will express in terms of~$U$ the generating function $P_{(\s,\lab^\star)}$ of pointed triangulations whose encoding mobile has normalized scheme $(\s,\lab^\star)$. We start by arbitrarily associating with every half-edge of~$\s$ an incident corner, and we do this in a bijective way. Let $(\fu,\lab)$ be a mobile encoding a triangulation under consideration. Observe that every corner of~$\s$ corresponds to a (possibly empty) tree of~$\fu$. We decompose~$\fu$ at the vertices corresponding to the vertices of~$\s$ without detaching the previous trees from the associated half-edge. We thus decompose $(\fu,\lab)$ into a collection of objects that are counted by some $C_{i,j}$'s, one such object per edge of~$\s$. Up to the location of the root, this gives an unequivocal decomposition of~$(\fu,\lab)$. As before, we need to count rooted mobiles twice and edge-rooted mobiles once, which amounts to counting nonrooted mobiles with a weight of double their number of edges. We let~$E$ be the edge set of~$\s$ and, for every edge $e\in E$, we denote by $\im(e)$, $\ip(e)\in\{1,\ldots,k\}$ the indices such that $\lab^\star_{\im(e)}\le\lab^\star_{\ip(e)}$ are the labels of the extremities of~$e$. We thus obtain
$$P_{(\s,\lab^\star)}=\frac{1}{2\,|E|}\,2t\frac{d}{dt}\left(\sum_\lab\prod_{e\in E} C_{\lab_{\im(e)},\lab_{\ip(e)}}\right),$$
where the sum is over all labelings whose normalization is~$\lab^\star$. The operator $2t\frac{d}{dt}$ counts the number of edges twice and the prefactor comes from the fact that there are $2|E|$ possible re-rootings of~$\s$.

We define the following sets:
\begin{align*}
E_=^\circ	&=\big\{e\in E\,:\,\lab^\star_{\im(e)}=\lab^\star_{\ip(e)}\in\Z\big\},\\
E_=^\oov	&=\Big\{e\in E\,:\,\lab^\star_{\im(e)}=\lab^\star_{\ip(e)}\in\Z+\frac12\Big\},\\
E_{\neq}	&=\big\{e\in E\,:\,\lab^\star_{\im(e)}\neq\lab^\star_{\ip(e)}\big\},\\
V^\oov	&=\Big\{v\in V(\s)\,:\,\lab^\star(v)\in\Z+\frac12\Big\}.
\end{align*}
Using~\eqref{Cij} and observing that $\sum_{e\in E} |\{\lab^\star_{\im(e)},\lab^\star_{\ip(e)}\}\cap(\Z+\frac12)|=3|V^\oov|$ (because every vertex in $V^\oov$ has degree~$3$), we may express the previous sum as follows:
$$\sum_\lab\prod_{e\in E} C_{\lab_{\im(e)},\lab_{\ip(e)}}
	=(B-1)^{|E_=^\circ|}\,B^{|E_=^\oov|}\, B^{|E_{\neq}|}\,N^{\frac32|V^\oov|}\sum_\lab\prod_{e\in E_{\neq}} U^{\lab_{\ip(e)}-\lab_{\im(e)}}.$$
Now, the sum in the right-hand side above is equal to
\begin{align*}
	\sum_\lab\prod_{e\in E_{\neq}} \prod_{j=\im(e)}^{\ip(e)-1}U^{\lab_{j+1}-\lab_{j}}
	&=\sum_{\delta_1,\ldots,\delta_{k-1}\ge 0}\,\prod_{e\in E_{\neq}} \prod_{j=\im(e)}^{\ip(e)-1}U^{\lab^\star_{j+1}-\lab^\star_{j}+\delta_j}\\
	&=\sum_{\delta_1,\ldots,\delta_{k-1}\ge 0}\, \prod_{j=1}^{k-1}U^{d(j)(\lab^\star_{j+1}-\lab^\star_{j}+\delta_j)}\\
	&=\prod_{j=1}^{k-1}\frac{U^{d(j)(\lab^\star_{j+1}-\lab^\star_{j})}}{1-U^{d(j)}}
\end{align*}
where we set $d(j)\de |\{e\in E_{\neq}\,:\,\im(e)\le j<\ip(e)\}|$ for $1\le j<k$.
Summing up and expressing everything in terms of~$U$ with the help of~\eqref{ZBU}, we obtain
\begin{equation}\label{Psl}
P_{(\s,\lab^\star)}=\frac{2^{|E_=^\circ|}}{|E|}\,3x\frac{d}{dx}\left(\frac{U^{|E_=^\circ|+\frac12|V^\oov|}\,(1+U)^{|E|-|E_=^\circ|-|V^\oov|}}{(1-U)^{|E|}}\prod_{j=1}^{k-1}\frac{U^{d(j)(\lab^\star_{j+1}-\lab^\star_{j})}}{1-U^{d(j)}}\right).
\end{equation}

\medskip
By the Euler characteristic formula, a triangulation of~$\cS$ with~$3n$ edges has $n+2-2h$ vertices. The generating functions of pointed triangulations and of triangulations of~$\cS$ are thus respectively
\begin{equation}\label{gfuntri}
\sum_{(\s,\lab^\star)\in\Sg} x^{2-2h} P_{(\s,\lab^\star)}\qquad\text{ and }\qquad \sum_{(\s,\lab^\star)\in\Sg} \int x^{1-2h}P_{(\s,\lab^\star)}\,dx\,.
\end{equation}

The leading contributions in the asymptotic formula will come from normalized schemes $(\s,\lab^\star)$ that maximize both the number of edges of~$\s$ and the cardinality of the range of~$\lab^\star$. It is not very hard to see that such normalized schemes are the ones where the scheme is cubic, that is, with vertices of degree~$3$ only, and the labeling function is injective. Moreover, such schemes all have $6h-3$ edges and $4h-2$ vertices. See \cite[Lemma~4.3]{ChDo15bij} for a proof of these facts. We denote by $\Sg^\star \subseteq\Sg$ the set of these normalized schemes. 

\begin{prop}\label{cS}
We suppose here that $h\ge 1$. The constant $c_\cS$ of Proposition~\ref{enum} is given by
$$c_\cS= \frac{2^{-(13h-9)/2}\,3^{-5(h-1)/2}}{(6h-3)\,\Gamma\big((5h-3)/2\big)}\Bigg(\sum_{(\s,\lab^\star)\in\Sg^\star}2^{-|V^\oov|}\prod_{j=1}^{4h-3}\frac{1}{d(j)}\Bigg).$$
\end{prop}

\begin{pre}[Proof of Proposition~\ref{enum} when $h\ge 1$ and of Proposition~\ref{cS}]
Let $(\s,\lab^\star)\in\Sg$ and let~$k$ denote the cardinality of the range of~$\lab^\star$ as above. The generating function $P_{(\s,\lab^\star)}$ admits a unique dominant singularity at $x=(12\sqrt{3})^{-1}$, which corresponds to a singularity of~$U$ and at the same time to a value where $U=1$. Moreover, at this singularity, $1-U$ admits the following Puiseux expansion:
$$1-U=2^{5/4}\, 3^{1/4} \big(1-12\sqrt{3}\,x\big)^{1/4}+\OO\big(\big(1-12\sqrt{3}\,x\big)^{1/2}\big)\,.$$
As a result,
$$[x^n] P_{(\s,\lab^\star)}\sim\frac{2^{-(|E|+4|V^\oov|+5k-5)/4}\,3^{-(|E|+k-5)/4}}{|E|\,\Gamma\big((|E|+k-1)/4\big)}\left(\prod_{j=1}^{k-1}\frac{1}{d(j)}\right) n^{(|E|+k-1)/4} \big(12\sqrt 3\big)^n\,.$$
The leading contributions are thus obtained for normalized schemes maximizing $|E|+k$ as we claimed. The number of triangulations of~$\cS$ with~$3n$ edges is thus asymptotically equivalent to
\begin{align*}
\frac1n \sum_{(\s,\lab^\star)\in\Sg} [x^n] P_{(\s,\lab^\star)}
	&\sim\frac{2^{-(13h-9)/2}\,3^{-5(h-1)/2}}{(6h-3)\,\Gamma\big((5h-3)/2\big)}\Bigg(\sum_{(\s,\lab^\star)\in\Sg^\star}2^{-|V^\oov|}\prod_{j=1}^{4h-3}\frac{1}{d(j)}\Bigg)\, n^{5(h-1)/2} \big(12\sqrt 3\big)^n
\end{align*}
as desired.
\end{pre}

For $h=1$, the expression of~$c_\cS$ becomes
$$c_\cS= \frac1{12}\Bigg(\sum_{(\s,\lab^\star)\in\Sg^\star}\frac1{2^{|V^\oov|}\,d(1)}\Bigg).$$
If~$\cS$ is the torus, there is only one possible scheme, consisting of~$3$ edges linking~$2$ vertices. In this case, $d(1)=3$ for every scheme and the possible labelings, starting from the root vertex, are $(0,1)$, $(0,\frac12)$, $(\frac12,1)$, $(\frac12,\frac32)$ and the symmetrical pairs. As a result,
$$c_\cS=\frac1{12}\frac23 \left(1+\frac12+\frac12+\frac14\right)=\frac18.$$

If~$\cS$ is the Klein bottle, there are two types of possible schemes obtained as follows. We consider a hexagon with sides denoted by $s_1$, $s_2$, \ldots, $s_6$ in clockwise order around it. We glue together~$s_1$ with~$s_4$ in an orientable way. The first type of schemes is then obtained by gluing in a nonorientable way~$s_2$ with~$s_6$ and~$s_3$ with~$s_5$. The second type of schemes is obtained by gluing in a nonorientable way~$s_2$ with~$s_3$ and~$s_5$ with~$s_6$. For the first type of scheme, $d(1)=3$ and for the second type of scheme, there are two loops so that $d(1)=1$. In both cases, the labels can be $\{0,1\}$, $\{0,\frac12\}$, $\{\frac12,1\}$ or $\{\frac12,\frac32\}$ with~$6$ possible rootings. All in all,
$$c_\cS=\frac1{12}\,6\left(1+\frac13\right) \left(1+\frac12+\frac12+\frac14\right)=\frac32.$$
This completes the values given in Table~\ref{tabcS}.

\begin{pre}[Proof of Proposition~\ref{gfun} when~$h=1$]
By~\eqref{Psl} and~\eqref{gfuntri}, the generating function of triangulations when $h=1$ is
$$\sum_{(\s,\lab^\star)\in\Sg}3\,\frac{2^{|E_=^\circ|}\,U^{|E_=^\circ|+\frac12|V^\oov|}\,(1+U)^{|E|-|E_=^\circ|-|V^\oov|}}{|E|\,(1-U)^{|E|}}\prod_{j=1}^{k-1}\frac{U^{d(j)(\lab^\star_{j+1}-\lab^\star_{j})}}{1-U^{d(j)}}.$$
Let us start with the torus. There are $11$ normalized schemes on the torus. There is one with one vertex labeled~$0$ and~$2$ edges: its contribution in the sum is ${6U^2}(1-U)^{-2}$.
The~$10$ remaining ones all have the same scheme with~$2$ vertices linked by~$3$ edges. The different labelings and the corresponding contributions are presented in Table~\ref{tabscheme}.

\begin{table}[ht!]%
\begin{center}
\begin{tabular}{cccccc}
$\big(\frac12,\frac12\big)$	& $\big(\frac12,\frac32\big)$, $\big(\frac32,\frac12\big)$	
	& $\big(0,\frac12\big)$, $\big(\frac12,0\big)$, $\big(1,\frac12\big)$, $\big(\frac12,1\big)$
										& $\big(0,0\big)$			& $\big(0,1\big)$, $\big(1,0\big)$	\\[2mm]
\hline\hline\\[-2mm]
$\dfrac{U(1+U)}{(1-U)^3}$	& $\dfrac{U^4(1+U)}{(1-U^3)(1-U)^3}$
	& $\dfrac{U^2(1+U)^2}{(1-U^3)(1-U)^3}$		& $\dfrac{8U^3}{(1-U)^3}$	& $\dfrac{U^3(1+U)^3}{(1-U^3)(1-U)^3}$
\end{tabular}
\end{center}
\caption{The contributions of all the normalized schemes on the torus except the one with only one vertex. The labeling is given as a pair whose first coordinate is the label of the root vertex and the second coordinate is the label of the other vertex.}
\label{tabscheme}
\end{table}

Summing these contributions, we obtain that the generating function of triangulations on the torus is equal to
$$\frac{U(U^2+10U+1)}{(1-U)^4}=\frac{\si(1-\si)}{2\,\big(1-6\si+6\si^2\big)^2}\,.$$
%
%

\medskip
On the Klein bottle, the normalized schemes are as follows. There is one with one vertex labeled~$0$ and two twisted edges and two with one vertex labeled~$0$, one twisted edge and one straight edge. These three normalized schemes all have contribution ${6U^2}(1-U)^{-2}$. The remaining normalized schemes all have a scheme of one of the two types presented before this proof. Moreover, for any such normalized scheme, nonrooted but with a prescribed root vertex, there are three possible rootings. The contributions of the normalized schemes of the first type are the same as for the torus, that is, the ones in Table~\ref{tabscheme}. For the normalized schemes of the second type, the contributions are presented in Table~\ref{tabscheme2}.

\begin{table}[ht!]%
\begin{center}
\begin{tabular}{cccccc}
$\big(\frac12,\frac12\big)$	& $\big(\frac12,\frac32\big)$, $\big(\frac32,\frac12\big)$	
	& $\big(0,\frac12\big)$, $\big(\frac12,0\big)$, $\big(1,\frac12\big)$, $\big(\frac12,1\big)$
										& $\big(0,0\big)$			& $\big(0,1\big)$, $\big(1,0\big)$	\\[2mm]
\hline\hline\\[-2mm]
$\dfrac{U(1+U)}{(1-U)^3}$	& $\dfrac{U^2(1+U)}{(1-U)^4}$
	& $\dfrac{2U^2(1+U)}{(1-U)^4}$		& $\dfrac{8U^3}{(1-U)^3}$	& $\dfrac{4U^3(1+U)}{(1-U)^4}$
\end{tabular}
\end{center}
\caption{The contributions of the normalized schemes of the second type on the Klein bottle. The labeling is given as a pair whose first coordinate is the label of the root vertex and the second coordinate is the label of the other vertex.}
\label{tabscheme2}
\end{table}

The generating function of triangulations on the Klein bottle is thus equal to
\begin{align*}
\frac{6U(13U^2+10U+1)}{(1-U)^4}&=\frac{3\si(1-\si)\big(7-30\si+30\si^2-6(1-2\si)\sqrt{1-6\si+6\si^2}\big)}{\big(1-6\si+6\si^2\big)^2}\,. \\[-4.5mm]
& \qedhere
\end{align*}
%
%
\end{pre}

\appendix

\section{Alternate orientation processes}\label{secquad}

We worked in this paper with the coherent orientation but, in fact, we may consider other orientation processes. Let us consider a triple $(\fu,\lab,\cO)$, where $(\fu,\lab)$ is a labeled unicellular mobile and~$\cO$ is the data of an orientation for each of its corner cycles. A corner is said to be \emph{$\cO$-oriented} if it is oriented as its corner cycle and the \emph{$\cO$-successor} $\suc_\cO(c)$ of a corner~$c$ with label $\lab(c)\ge 2$ is the first subsequent corner with label strictly smaller than~$\lab(c)$ in the order given by the $\cO$-orientation of~$c$. Finally, the triple $(\fu,\lab,\cO)$ is a \emph{well-labeled loop-oriented unicellular mobile} if, for every corner~$c$ with label $\lab(c)\ge 2$, we have $\lab(\suc_\cO(c))=\lab(c)-1$.

On the other side, a \emph{loop-oriented pointed bipartite map} is a triple $(\m,v^\ooo,\cO)$ consisting of a pointed bipartite map $(\m,v^\ooo)$ with the additional data of an orientation for each of its level loops that are not maximal level loops. The sets of well-labeled loop-oriented unicellular mobiles and loop-oriented pointed bipartite maps are in $2$-to-$1$ correspondence through the following bijections.

From a well-labeled loop-oriented unicellular mobile $(\fu,\lab,\cO)$ and a parameter $\eps\in\{+,-\}$, we define the pointed bipartite map $(\m,v^\ooo)$ by the construction of Section~\ref{secum}, replacing successors by $\cO$-successors. The facts that this construction is well defined and outputs a pointed bipartite map are proved exactly as with the coherent orientation. As noticed during the proof of Theorem~\ref{bijthm}, there is a one-to-one correspondence between the corner cycles of~$(\fu,\lab)$ and the level loops of~$(\m,v^\ooo)$ that are not maximal level loops. We orient the level loops of~$(\m,v^\ooo)$ that are not maximal level loops in accordance with the corresponding corner cycles, as in Figure~\ref{ccll}.

From a loop-oriented pointed bipartite map $(\m,v^\ooo,\cO)$, we define the unicellular mobile $(\fu,\lab)$ by the construction of Section~\ref{secmu}, using the $\cO$-orientation of its level loops instead of the process of Section~\ref{secori}. The fact that the maximal level loops are not oriented does not bear any effect as both possible orientations for a maximal level loop set the same stops for the loop. From the proof of Theorem~\ref{bijthm}, the corner cycles of~$(\fu,\lab)$ correspond to the level loops of~$(\m,v^\ooo)$ that are not maximal level loops: we orient the corner cycles in accordance with the corresponding level loops.

The proof that the above mappings are bijections is easily obtained from the arguments of Section~\ref{secpfinv}. As a result, using another orientation process for the level loops of a pointed bipartite map may provide a bijection with a possibly different set of encoding objects. In the case of pointed bipartite quadrangulations (see Corollary~\ref{bijcor} in the case $\alpha=(2,2,\ldots,2)$), the situation is particularly nice because we obtain well-labeled unicellular mobiles with green vertices of degree~$2$ (which are in direct bijection with well-labeled unicellular maps, defined in Section~\ref{secencquad}). As explained at the end of Section~\ref{seccoh}, the condition of being well labeled in this case does not depend on the way the corners are oriented. Consequently, any orientation process of the level loops of the bipartite quadrangulation yields a bijection, as long as this orientation process is tractable on the encoding unicellular maps.

If we use the restrictions of~$\Phi$ and~$\Psi$ defined in Sections~\ref{secmu} and~\ref{secum} in the case of bipartite quadrangulations, some simplifications appear. There are no trivial corner runs, as for a white oriented corner~$\vc$, we have $\varphi\circ\si\circ\varphi(\vc)=\si(\vc)$ so that the labels of~$\vc$ and of $\varphi\circ\si\circ\varphi(\vc)$ are equal. The level loop also has a simpler structure: when constructing a level~$i$-loop, if we cross an \edg{$i$}{$i$+$1$} edge and move along an \edg{$i$}{$i$-$1$} edge, it means that we are visiting a face whose labels are $i-1$, $i$, $i+1$, $i$, so that the next edge we encounter is an \edg{$i$-$1$}{$i$} edge and we move along it. As a result, the corner labels along the level loop are several~$i$'s, then one $i-1$, then several~$i$'s and so on. The maximal level loops visit four corners, labeled $i$, $i-1$, $i$, $i-1$ or $i$, $i-1$, $i-2$, $i-1$. 

\bigskip

Instead of describing the restrictions of~$\Phi$ and~$\Psi$, let us give alternate orientation processes for level loops and corners. This is the option we proposed in~\cite{bettinelli15fpsac}.

\subsection{From pointed bipartite quadrangulations to well-labeled unicellular maps}\label{secqu}

Let us consider a pointed bipartite quadrangulation $(\q,v^\ooo)$. We perform the steps of Section~\ref{secmu}, except that, in Step~\ref{step2}, we use the following orientation process of the level loops instead of that of Section~\ref{secori}.

\begin{figure}[ht!]
		\psfrag{v}[][]{$v^\ooo$}
		\psfrag{0}[][][.8]{$0$}
		\psfrag{1}[][][.8]{$1$}
		\psfrag{2}[][][.8]{$2$}
		\psfrag{3}[][][.8]{$3$}
		\psfrag{4}[][][.8]{$4$}
		\psfrag{a}[][][.8]{$a$}
		\psfrag{b}[][][.8]{$b$}
		\psfrag{c}[][][.8]{$c$}
		\psfrag{d}[][][.8]{$d$}
		\psfrag{e}[][][.8]{$e$}
		\psfrag{f}[][][.8]{$f$}
		\psfrag{g}[][][.8]{$g$}
		\psfrag{h}[][][.8]{$h$}
		\psfrag{i}[][][.8]{$i$}
		\psfrag{j}[][][.8]{$j$}
		\psfrag{k}[][][.8]{$k$}
		\psfrag{l}[][][.8]{$l$}
		\psfrag{m}[][][.8]{$m$}
	\centering\includegraphics[width=.95\linewidth]{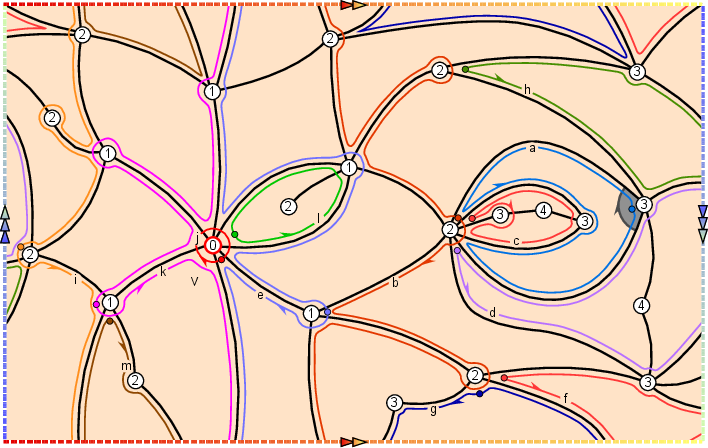}
	\caption{Level loops of a pointed bipartite quadrangulation. We did not represent the 18 maximal level loops of simple faces. The letters indicate the rank of the loops in the order (not considering maximal level loops of simple faces), the half arrowheads their orientation and the dots their origin. Loop~$j$ is at level~$0$, loops $e$, $k$ and~$l$ are at level~$1$, loops~$b$, $i$ and~$m$ are at level~$2$ and the remaining loops are at level~$3$.}
	\label{loopsq}
\end{figure}

We use the root~$\vrho$ of the quadrangulation in order to orient, give an origin to and order all the level loops. The first loops in our order are the ones that visit the root corner, ordered by decreasing level (there may be~$1$ or~$2$ such loops). We orient them in the orientation of~$\vrho$ and set their origin at the location of the root corner. Next, we start from the origin of the first loop and travel on it, according to its orientation. Every time we encounter a new loop (in the sense that we visit a corner that is also visited by another loop that has not yet been oriented), we declare it to be the next one in the order we are creating. We also set its origin at the location at which we are and we orient it in the orientation induced by the loop on which we are traveling. When we arrive back at the origin of our loop, we move to the next one in the order and iterate the process until every loop has been oriented; see Figure~\ref{loopsq}.

This operation terminates as every loop at level $i\ge 1$ visits at least a corner visited by a loop at level $i-1$. Indeed, a loop issued from an oriented corner~$\vc$ visits the oriented corner~$\pr(\vc)$, which is visited by the loop issued from it.

\begin{rem}
If the surface we consider is orientable, all the loops are merely oriented according to the orientation of the surface induced by the orientation of the root.
\end{rem}

Observe that, in a so-called \emph{simple face}, that is, a face with corner labels $i$, $i+1$, $i+2$, $i+1$, there are three level loops:
\begin{itemize}
	\item one at level~$i$ that visits the corner \corn{$i$+$1$}{$i$}{$i$+$1$};
	\item one at level $i+1$ that visits the corners \corn{$i$+$2$}{$i$+$1$}{$i$}, \corn{$i$+$1$}{$i$}{$i$+$1$} and \corn{$i$}{$i$+$1$}{$i$+$2$};
	\item and one maximal level loop that visits all the corners.
\end{itemize}
Consequently, not considering the maximal level loops of the simple faces in the above process does not change the orientations, origins, and relative order of the other level loops. We thus discard these level loops, except possibly the one that visits the root corner, if it is the only level loop that visits the root corner. 

\begin{figure}[ht!]
		\psfrag{0}[][][.8]{$0$}
		\psfrag{1}[][][.8]{$1$}
		\psfrag{2}[][][.8]{$2$}
		\psfrag{3}[][][.8]{$3$}
		\psfrag{4}[][][.8]{$4$}
		\psfrag{v}[][]{$v^\ooo$}
	\centering\includegraphics[width=.95\linewidth]{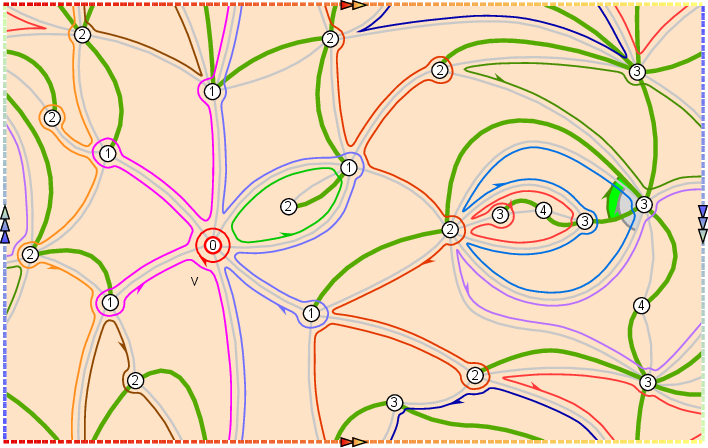}
	\caption{The bijection, from a pointed bipartite quadrangulation to a well-labeled unicellular map.}
	\label{bij_qu}
\end{figure}

Finally, we remove the green vertices, which are all of degree~$2$. Figure~\ref{bij_qu} shows an example of the construction.

\subsection{From well-labeled unicellular maps to pointed bipartite quadrangulations}\label{secuq}

The inverse mapping constructs a pointed bipartite quadrangulation from a well-labeled unicellular map $(\fu,\lab)$ and a parameter $\eps\in\{+,-\}$ as follows. 

\begin{figure}[ht!]
		\psfrag{0}[][][.8]{$0$}
		\psfrag{1}[][][.8]{$1$}
		\psfrag{2}[][][.8]{$2$}
		\psfrag{3}[][][.8]{$3$}
		\psfrag{4}[][][.8]{$4$}
		\psfrag{v}[][]{$v^\ooo$}
	\centering\includegraphics[width=.95\linewidth]{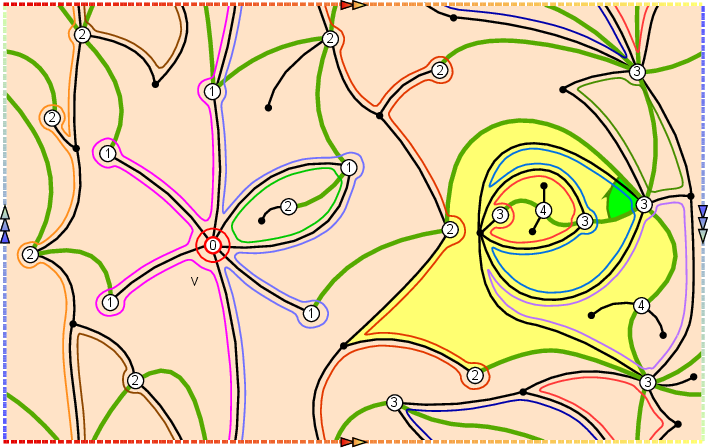}
	\caption{The sectors and level loops. The sector at level~$2$ containing the root of~$\fu$ is highlighted. See Figure~\ref{loopsuq_p} for a polygon representation of the same figure.}
	\label{loopsuq}
\end{figure}

\begin{figure}[ht!]
\floatbox[{\capbeside\thisfloatsetup{capbesideposition={left,center},floatwidth=sidefil}}]{figure}[.65\textwidth]
{\caption{Polygon representation of Figure~\ref{loopsuq}. We represented the identifications of the loops at level~$1$ on the paired sides by dashed lines (regardless of the orientability of the pairings). The root of the future quadrangulation is represented in gray.}}
{		\psfrag{0}[][][.7]{$0$}
		\psfrag{1}[][][.7]{$1$}
		\psfrag{2}[][][.7]{$2$}
		\psfrag{3}[][][.7]{$3$}
		\psfrag{4}[][][.7]{$4$}
	\centering\includegraphics[width=.95\linewidth]{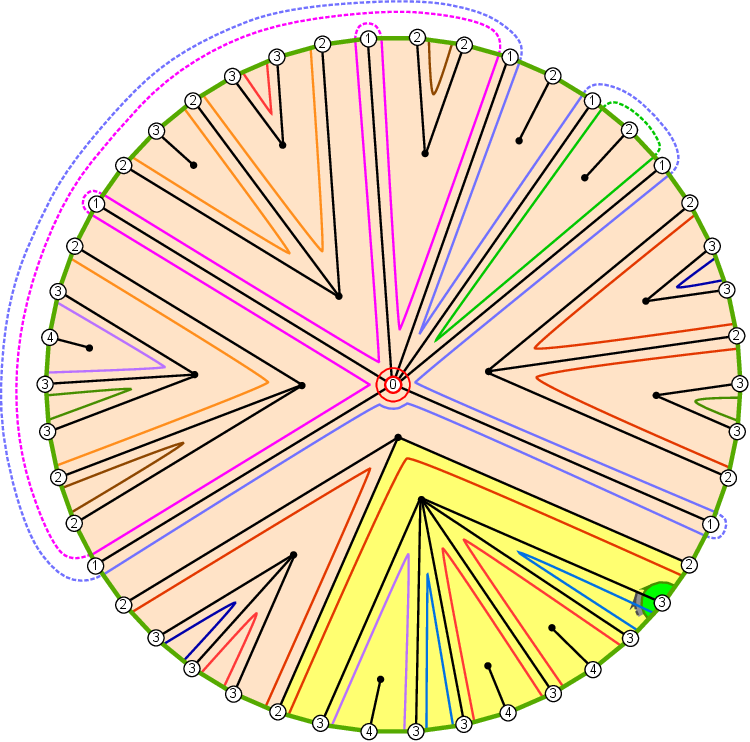}
	\label{loopsuq_p}
}
\end{figure}

\paragraph{Step 1. Defining the sectors [Figures~\ref{loopsuq} and~\ref{loopsuq_p}].}
We first add inside the unique face of~$\fu$ a new vertex~$v^\ooo$ and assign to it the label $\lab(v^\ooo)\de 0$. We connect all the corners with label~$1$ to this vertex~$v^\ooo$ in a noncrossing fashion. We thus create a certain number of \emph{sectors}, defined as the connected components of the complement of the newly added edges and the original edges of~$\fu$. The set of corners of~$\fu$ belonging to a sector is called the \emph{arc} of the sector and the remaining corner of the sector is called its \emph{inner corner}.

Inside each sector whose arc contains at least three corners, we add a temporary vertex to which we connect all the corners with label~$2$ (note that the first and last corners in such a sector are labeled~$1$ by definition, and the second and second to last are necessarily labeled~$2$). The added edges delimit new sectors, defined as the connected components that do not contain the inner corner of the original sector. We iterate the construction inside each created sector. The \emph{level} of a sector is the minimal label of the corners on its arc (attained only once at each extremity of the arc, by definition).

\paragraph{Step 2. Constructing the level loops [Figures~\ref{loopsuq} and~\ref{loopsuq_p}].}
A level loop is constructed as follows. We start from the inner corner of a sector and move along one of the two edges linking it to a corner of~$\fu$. When we reach an edge of~$\fu$, we cross it and continue to move along the new edge we meet, toward the inner corner of the sector we visit. When we reach the inner corner, we move along the side of the other edge linking it to a corner of~$\fu$, and we iterate the process until we close the loop. The \emph{level} of a loop is the common level of the sectors it visits. We also add around~$v^\ooo$ a single loop at level~$0$.

Before orienting the loops, we set the root~$\vrho$ of the future quadrangulation by the rule depicted in Figure~\ref{rootm}. In the eventuality that no loop visits~$\vrho$, we add a loop that circles around the middle of the green edge preceding the root of~$\fu$. 

\begin{figure}[ht!]
		\psfrag{2}[][][.8]{$2$}
		\psfrag{3}[][][.8]{$3$}
		\psfrag{4}[][][.8]{$4$}		
	\centering\includegraphics[width=.95\linewidth]{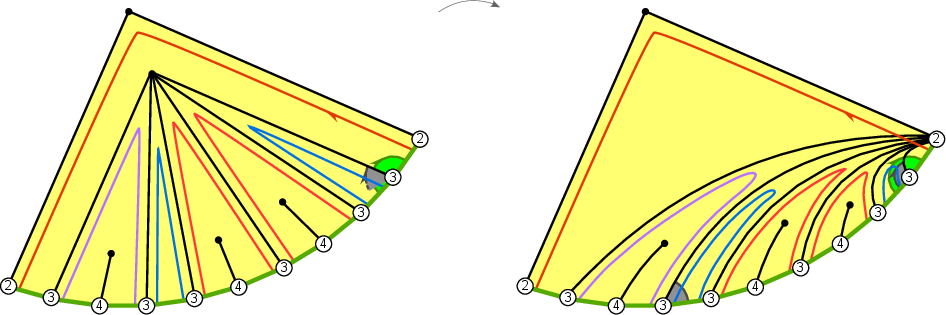}
	\caption{Identification of the temporary vertex of the highlighted sector of Figures~\ref{loopsuq} and~\ref{loopsuq_p}. When the dark red level~$2$-loop is oriented (from right to left in the figure), we identify the temporary vertex of the level~$3$-subsectors with the extremity of the arc first visited by the portion of loop (the one at the right in the figure).}
	\label{identify}
\end{figure}

\paragraph{Step 3. Orienting the level loops and identifying the temporary vertices.}
We now orient, give an origin to and order the level loops. At the same time, we also identify the temporary vertices with the appropriate vertices of~$\fu$. More precisely, every time we orient a loop, inside each sector it visits, we identify the temporary vertex of the subsectors at next level with the first vertex of the sector that is visited by the portion of the level loop inside the sector; see Figure~\ref{identify}.

\begin{figure}[ht!]
		\psfrag{v}[][]{$v^\ooo$}
		\psfrag{a}[][][.8]{$a$}
		\psfrag{b}[][][.8]{$b$}
		\psfrag{c}[][][.8]{$c$}
		\psfrag{d}[][][.8]{$d$}
		\psfrag{e}[][][.8]{$e$}
		\psfrag{f}[][][.8]{$f$}
		\psfrag{g}[][][.8]{$g$}
		\psfrag{h}[][][.8]{$h$}
		\psfrag{i}[][][.8]{$i$}
		\psfrag{j}[][][.8]{$j$}
		\psfrag{k}[][][.8]{$k$}
		\psfrag{l}[][][.8]{$l$}
		\psfrag{m}[][][.8]{$m$}
		\psfrag{0}[][][.8]{$0$}
		\psfrag{1}[][][.8]{$1$}
		\psfrag{2}[][][.8]{$2$}
		\psfrag{3}[][][.8]{$3$}
		\psfrag{4}[][][.8]{$4$}
	\centering\includegraphics[width=.95\linewidth]{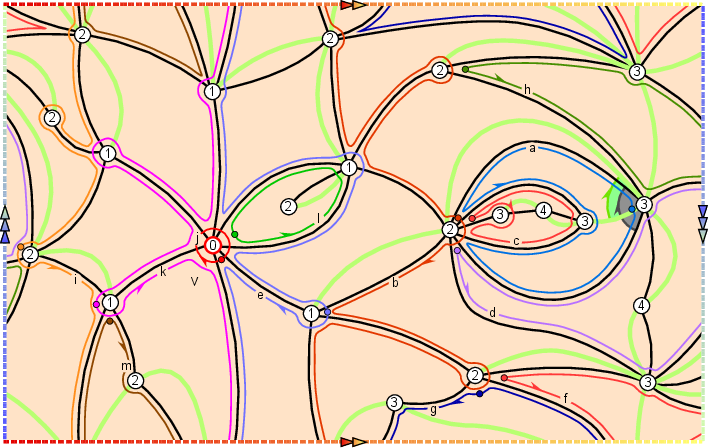}
	\caption{The bijection, from a well-labeled unicellular map to a pointed bipartite quadrangulation. Note that everything was done in such a way that the level loops correspond through the bijection. See Figure~\ref{bij_uq_p} for a polygon representation of the same figure.}
	\label{bij_uq}
\end{figure}

\begin{figure}[ht!]
\floatbox[{\capbeside\thisfloatsetup{capbesideposition={right,center},floatwidth=sidefil}}]{figure}[.65\textwidth]
{\caption{Polygon representation of Figure~\ref{bij_uq}. The bijection, from a well-labeled unicellular map to a pointed bipartite quadrangulation.}}
{		
		\psfrag{a}[][][.7]{$a$}
		\psfrag{b}[][][.7]{$b$}
		\psfrag{c}[][][.7]{$c$}
		\psfrag{d}[][][.7]{$d$}
		\psfrag{e}[][][.7]{$e$}
		\psfrag{f}[][][.7]{$f$}
		\psfrag{g}[][][.7]{$g$}
		\psfrag{h}[][][.7]{$h$}
		\psfrag{i}[][][.7]{$i$}
		\psfrag{j}[][][.7]{$j$}
		\psfrag{k}[][][.7]{$k$}
		\psfrag{l}[][][.7]{$l$}
		\psfrag{m}[][][.7]{$m$}
		\psfrag{0}[][][.7]{$0$}
		\psfrag{1}[][][.7]{$1$}
		\psfrag{2}[][][.7]{$2$}
		\psfrag{3}[][][.7]{$3$}
		\psfrag{4}[][][.7]{$4$}
	\centering\includegraphics[width=.95\linewidth]{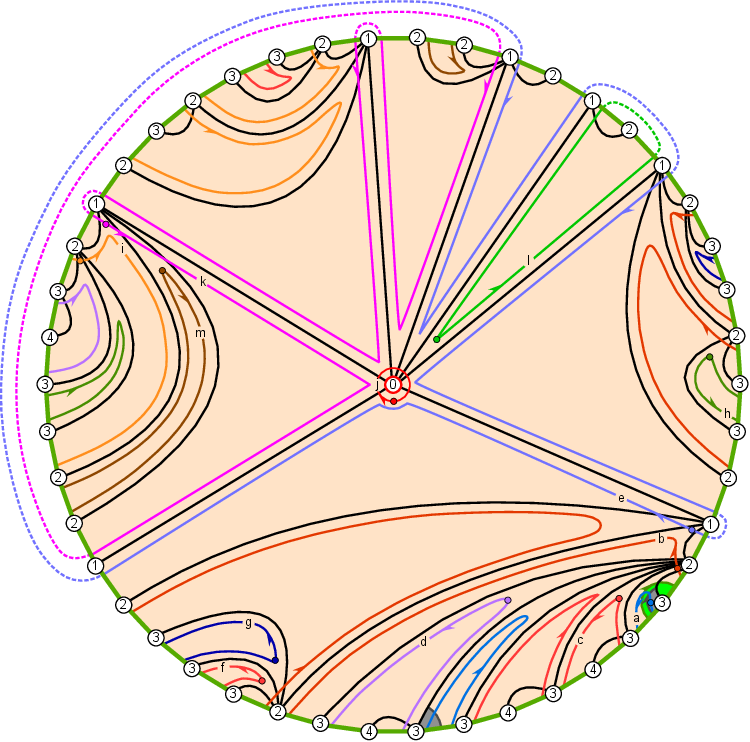}
	\label{bij_uq_p}
}
\end{figure}

The first loop in our order is the one that visits~$\vrho$. We set its origin at the location of~$\vrho$ and orient it in the orientation of~$\vrho$. We do the aforementioned identifications. We then start from the origin of the first loop and travel on it, according to its orientation. Every time we visit a temporary vertex that has not been identified, we orient the loop visiting the sector at previous level in the same orientation as the orientation of the loop on which we are traveling, and do the resulting identifications of temporary vertices. After the identifications of temporary vertices, the corner we are visiting is also visited by the loop we just oriented. We declare this loop to be the next one and we set its origin at the location we are. Similarly, every time we visit an inner corner and discover a new loop, we declare it to be the next one, we set its origin at the location we are and we orient it in the orientation induced by the loop on which we are traveling. When we arrive back at the origin of our loop, we move to the next one in our order and iterate the process until every loop has been oriented; see Figures~\ref{bij_uq} and~\ref{bij_uq_p}.

\paragraph{Step 4. Discarding the green edges.}
The embedded graph~$\q$ made of the added edges and rooted at the future root is a bipartite quadrangulation. If $\eps=-$, then the output of the construction is $(\q,v^\ooo)$, and if $\eps=+$, then the output is the root flipped version $(\bar\q,v^\ooo)$.

\subsection{Conclusion}

\begin{thm}\label{thmquad}
The mappings of Appendix~\ref{secquad} are inverse bijections.
\end{thm}

\begin{pre}
Let us consider a pointed bipartite quadrangulation $(\q,v^\ooo)$ and orient its level loops by the process of Section~\ref{secqu}. Let also~$(\fu,\lab)$ be the well-labeled unicellular map obtained from the construction of Section~\ref{secqu}. Then the well-labeled loop-oriented unicellular mobile corresponding to $(\q,v^\ooo)$ through the bijection described at the beginning of Appendix~\ref{secquad} is the labeled unicellular mobile with green vertices of degree~$2$ obtained from $(\fu,\lab)$ by adding a green vertex on every edge, and whose corner cycles are oriented in accordance with the corresponding level loops oriented by the process of Section~\ref{secuq}.

With our orientation processes, we thus created a $1$-to-$1$ correspondence between pointed bipartite quadrangulations and some loop-oriented pointed bipartite quadrangulations, as well as a $1$-to-$1$ correspondence between well-labeled unicellular maps and some well-labeled loop-oriented unicellular mobile. Moreover, the two subsets of loop-oriented pointed bipartite quadrangulations and well-labeled loop-oriented unicellular mobiles are in $1$-to-$1$ correspondence through the bijection of the beginning of Appendix~\ref{secquad}.
\end{pre}

\paragraph{Other choices of orientations, other bijections.}
In the previous construction, we had the opportunity to choose the orientation of every loop and we decided to use, by default, the orientation induced by the root or by the loop we were visiting. We can modify this rule as follows. We fix beforehand a sequence of rules ``$+$'' or ``$-$'' and, every time we orient a new loop, we reverse its orientation if and only if the current rule is ``$-$.''

For every choice of sequence, we obtain a different pair of mappings, which are inverse one from another. In fact, any loop orientation process that is tractable provides a bijection.

\section[Representation of unicellular maps as polygons]{Representation of unicellular maps as polygons with paired sides}\label{annpol}

\begin{figure}[ht!]
		\psfrag{a}[][][.8]{$a$}
		\psfrag{b}[][][.8]{$b$}
		\psfrag{c}[][][.8]{$c$}
		\psfrag{d}[][][.8]{$d$}
		\psfrag{e}[][][.8]{$e$}
		\psfrag{f}[][][.8]{$f$}
		\psfrag{g}[][][.8]{$g$}
		\psfrag{h}[][][.8]{$h$}
		\psfrag{i}[][][.8]{$i$}
		\psfrag{j}[][][.8]{$j$}
		\psfrag{k}[][][.8]{$k$}
		\psfrag{l}[][][.8]{$l$}
		\psfrag{m}[][][.8]{$m$}
		\psfrag{n}[][][.8]{$n$}
		\psfrag{o}[][][.8]{$o$}
		\psfrag{p}[][][.8]{$p$}
	\centering\includegraphics[width=12cm]{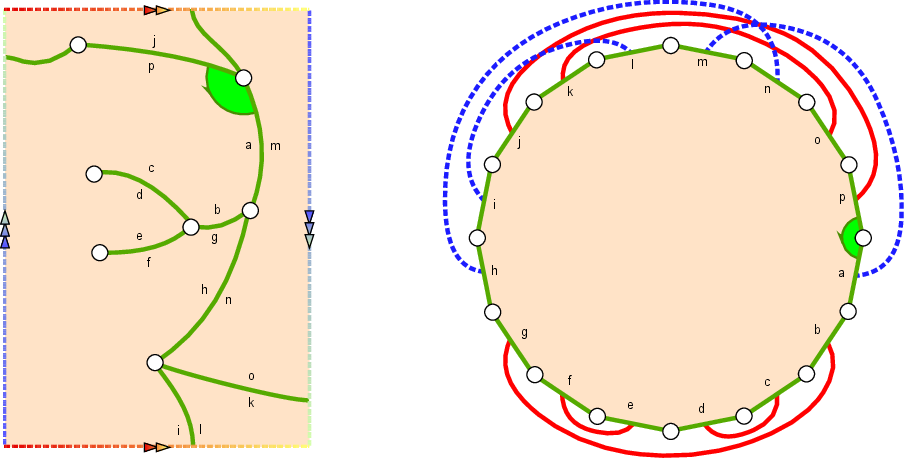}
	\caption{Representation of a unicellular map as a polygon with paired sides. The red straight lines correspond to orientable pairings and the blue dashed lines correspond to nonorientable pairings. The letters are positioned according to the root-induced contour order.}
	\label{polygon}
\end{figure}

It is sometimes convenient to represent a unicellular map~$\fu$ as a polygon with twice as many sides as the number of edges of~$\fu$. In this representation, the sides of the polygon are paired either in an orientable way or in a nonorientable way (see Figure~\ref{polygon}).

In this appendix, we use this representation for Figures~\ref{bij_um}, \ref{loopsuq} and~\ref{bij_uq}.


\begin{figure}[ht!]
\floatbox[{\capbeside\thisfloatsetup{capbesideposition={right,center},floatwidth=sidefil}}]{figure}[.65\textwidth]
{\caption{The bijection, from a well-labeled unicellular mobile to a pointed bipartite map (polygon representation of Figure~\ref{bij_um}). The oriented corner cycles and the geodesic orientation are represented. For better visibility, the side pairings are not represented. However, most of the identifications of the corner cycles on the paired sides are represented by dashed lines (regardless of the orientability of the pairings). Note that, in this example, there is only one corner whose coherent orientation differs from its root-induced orientation: it is the label~$4$-corner on the bottom left, in the purple corner cycle. For this reason, it is the only corner whose coherent orientation is represented. The half-arrowheads on the corner cycles are located at the positions of their first white or green corners in the root-induced contour order.}}
{		\psfrag{0}[][][.7]{$0$}
		\psfrag{1}[][][.7]{$1$}
		\psfrag{2}[][][.7]{$2$}
		\psfrag{3}[][][.7]{$3$}
		\psfrag{4}[][][.7]{$4$}
	\centering\includegraphics[width=.95\linewidth]{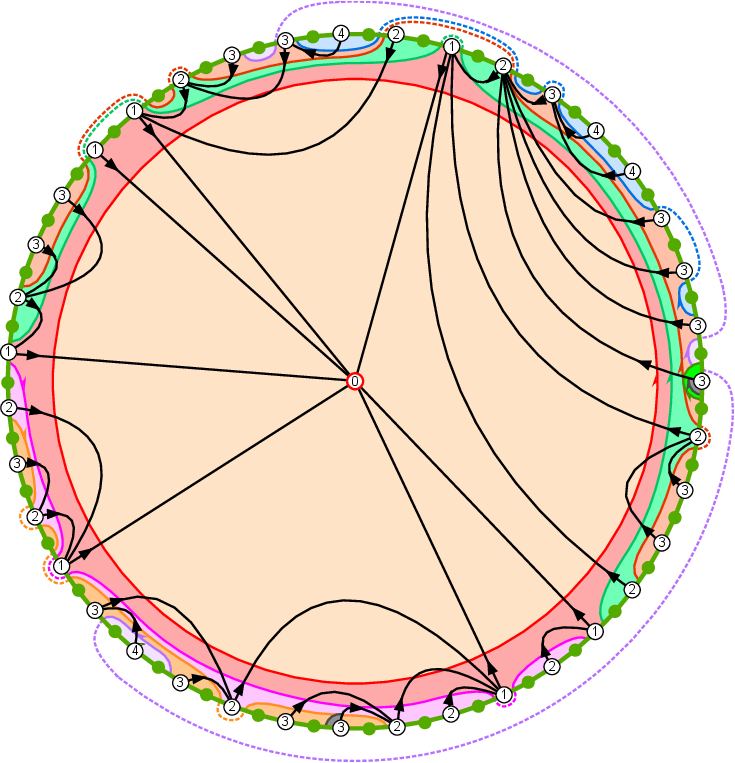}
	\label{bij_um_p}
}
\end{figure}


\clearpage

\phantomsection
\addcontentsline{toc}{section}{Often used notation}
\label{secnot}
\printnomenclature[15mm]

\bibliographystyle{alpha}
\bibliography{main}

\end{document}